\documentclass{amsart}

\pdfoutput=1
\usepackage{style}

\begin{document}

\title[]{Efficient computation of bifurcation diagrams with a deflated approach to reduced basis spectral element method}

\author[]{Moreno Pintore$^1$ \and Federico Pichi$^2$ \and Martin Hess$^2$ \and Gianluigi Rozza$^2$ \and Claudio Canuto$^1$}
\address{$^1$ Politecnico di Torino, Department of Mathematical Sciences ``G.L. Lagrange'', Corso Duca degli Abruzzi 24,  10129 Torino, Italy. Email: moreno.pintore@polito.it, claudio.canuto@polito.it\vspace{0.3cm}}
\address{$^2$ SISSA, International School for Advanced Studies, Mathematics Area, mathLab,
Via Bonomea 265, 34136 Trieste, Italy. Email: fpichi@sissa.it, mhess@sissa.it, gianluigi.rozza@sissa.it}

\maketitle

\begin{abstract}
The majority of the most common physical phenomena can be described using partial differential equations (PDEs). However, they are very often characterized by strong nonlinearities. Such features lead to the coexistence of multiple solutions studied by the bifurcation theory. Unfortunately, in practical scenarios, one has to exploit numerical methods to compute the solutions of systems of PDEs, even if the classical techniques are usually able to compute only a single solution for any value of a parameter when more branches exist. In this work we implemented an elaborated deflated continuation method, that relies on the spectral element method (SEM) and on the reduced basis (RB) one, to efficiently compute bifurcation diagrams with more parameters and more bifurcation points. The deflated continuation method can be obtained combining the classical continuation method and the deflation one: the former is used to entirely track each known branch of the diagram, while the latter is exploited to discover the new ones. Finally, when more than one parameter is considered, the efficiency of the computation is ensured by the fact that the diagrams can be computed during the online phase while, during the offline one, one only has to compute one-dimensional diagrams. In this work, after a more detailed description of the method, we will show the results that can be obtained using it to compute a bifurcation diagram associated with a problem governed by the Navier-Stokes equations.
\end{abstract}

\section{Introduction and motivation}
\begin{figure}
\centering
\includegraphics[width=139mm, angle=0, keepaspectratio]{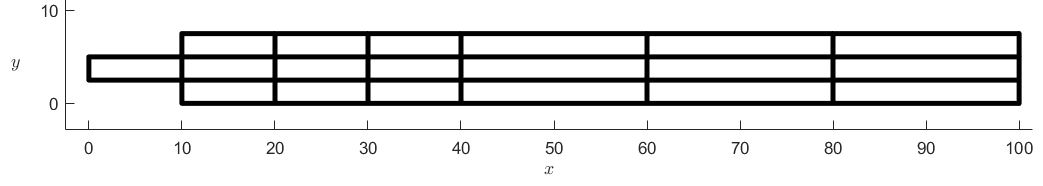}
\caption{Domain $\Omega$ and mesh used in the offline phase to compute the full order solutions}
\label{mesh}
\end{figure}

Usually, when one wants to numerically compute a bifurcation diagram, one has to combine many numerical methods in order to obtain it.  In fact, a suitable discretization method coupled with a nonlinear solver is required to compute any solution of the nonlinear problem of interest, while at least a continuation method and an additional technique to discover new branches are necessary to generate an entire diagram \cite{seydel}. However, since numerous solutions have to be computed and the involved system has to be solved several times, the computational cost of the task may be prohibitive in practical scenarios.  Note that the computational cost is often further increased because the unknown branches are frequently sought simply initializing the iterative solver with different guesses, trying to converge to new solutions. In order to decrease such a cost, we decided to implement a technique based  on an efficient combination of four different methods.\\
Firstly, we rely on the reduced basis method \cite{rozza_rom}. This is important because, after a very expensive offline phase, the computation of a solution in the online one is, in a repetitive computational environment, very efficient. In fact, the former is used to generate a low-dimensional space $V^{rb}$ defined as a combination of some of the most important solutions obtained during such a phase (these solutions are called snapshots or full order solutions). Subsequently, during the online phase, the solutions are sought in $V^{rb}$ and the affine decomposition of the operators can be exploited to further speed up the computation \cite{supremizer}. This way, it is possible to efficiently compute any solution and to discretize the entire bifurcation diagram only in the online phase, significantly reducing its associated computational cost.\\
Secondly, we decided to use the spectral element method (SEM) \cite{sem} to compute the snapshots in the offline phase. This is important because numerous solutions have to be computed to obtain a reduced space able to capture all the branches in the online phase. The offline phase is thus very expensive, but its cost is reduced exploiting the fact that the solutions obtained with the SEM are characterized by a lower number of degrees of freedom when compared to their counterpart computed with the standard finite element method (FEM). Moreover, the computational cost of the offline phase can be further decreased using the static condensation method (also known as Schur complement \cite{golub}) to efficiently solve the resulting linear system \cite{nektar_online}. \\
Finally, we implemented an elaborated deflated continuation method to compute the snapshots in the offline phase and the complete bifurcation diagram in the online one. Such a technique is composed by two different parts: the continuation method \cite{dijkstra} and the deflation one \cite{classic_deflation}. The former is used to follow a known branch of the diagram, while we used the latter to compute the first solutions on unknown branches. Note that the first solutions on unknown branches are the ones that cannot be computed with the continuation method because they do not belong to the branch of any of the previously computed solutions. Moreover, since they can be used to compute other solutions on their branches via the continuation method, they are the first ones belonging to such branches. The main idea behind the continuation method is to exploit the iterative solver and the continuity of each branch to obtain a solution very similar to the previous one and, therefore, belonging to the same branch of the latter. On the other hand, the deflation prevents the iterative solver from converging to known solutions and, this way, if it converges, it will converge to yet unknown fields.\\
In this work we only focus on the incompressible Navier-Stokes equations, although the SEM, the RB, the continuation and the deflation methods are numerical techniques that are not related to a particular class of equations and can thus be used in different frameworks. For instance, in \cite{pichi_schro} and in \cite{pichi_bif}, the authors analyzed the effectiveness of similar techniques (with the SEM substituted by the FEM and without the deflation) to compute bifurcation diagrams for the Gross-Pitaevskii and the {V}on {K}\'arm\'an equations.
In particular, the present work is strongly related to \cite{local_rom}, \cite{martin_geom}, \cite{martin_curve_walls} and \cite{martin0}. In fact, a possible future application of this work could be the mitral valve regurgitation \cite{cardio}. This is a cardiac disease characterized by an inverse blood flux from the left ventricle to the left atrium. The main tool to detect and analyze such a disease is echocardiography. However, when the blood flow undergoes the Coanda effect \cite{tritton}, it is very complex to quantify the flow rate. Therefore, it can be useful to use the direct simulation of the flow to properly analyze the images obtained via echocardiography.\\
We thus considered a two-dimensional channel with a narrow inlet, that represented, in a very simplified way, the mitral valve and the left atrium. Such a domain $\Omega$ is shown, with the mesh used in the offline phase, in figure \ref{mesh}.
Here the left vertical wall is the inlet, the right one is the outlet and the remaining ones represent the heart walls. It can easily be noted that the mesh is very coarse, in fact only 19 elements are involved. However, thanks to the SEM, it is possible to compute very accurate solutions even with meshes similar to this one. Indeed, it relies on high-order ansatz functions inside each element that allow an exponential decay of the error \cite{hp_conv}. On the contrary, if one had computed the same solutions with the FEM, a much finer mesh would have been required because of the algebraic convergence that characterizes such a method \cite{fem_error_b}. Moreover, since the convergence is much faster, it is possible to reach the same level of accuracy with a significantly lower number of degrees of freedom, thus further increasing the efficiency of the method, both in the assembly of the matrices and when the associated systems is solved.\\
It is important to remark that the results presented in section \ref{results_chapter} could be obtained using the deflated continuation method without the reduced basis approach and with an arbitrary discretization technique. However, the computational cost of such a process would be prohibitive because numerous solutions have to be computed to build a bifurcation diagram. In fact, it is important to highlight that the computational cost of the process exponentially increases with the parameter space dimension. Exploiting the described approach, instead, its cost is significantly reduced when several parameters are involved. In fact it is possible to compute few one-dimensional bifurcation diagrams during the offline phase and to reconstruct the other dimensions only in the online one, when each solution can be computed much more efficiently.\\
Finally, we highlight that there exist other techniques to compute bifurcation diagrams based on continuation methods (see \cite{dijkstra} or \cite{keller}) or expensive direct techniques as the branching system method \cite{seydel}. The advantage of the proposed technique over the previous ones is that it allows to efficiently and automatically detect new branches to compute an entire diagram. Moreover, we show, for the first time, that such a technique is stable also when combined with the RB method. Note that different approaches, mainly based on the analysis of the eigenvalues and the eigenvectors of the linearized equations, can be used exploiting an offline-online splitting as described, for instance, in \cite{pichi_schro} and \cite{pichi_bif}. Such papers, though, are more related to the stability analysis of the found solutions, while our approach aims at obtaining an entire diagram. Therefore these two techniques could be mixed to efficiently compute an entire bifurcation diagram while checking its stability properties.\\
The paper is organized as follows: in section \ref{prob_form} we will focus on the formulation of the problem of interest, deriving the weak formulation from the strong equations and presenting the two linearization techniques \cite{burger_notes} that will be used. Then, in section \ref{sem_chapter}, the SEM will be described, with a particular focus on the static condensation method. The latter is a technique that can be used to significantly speed up the computation of the solution of the involved linear system exploiting the fact that two kinds of modes are present. In section \ref{rb_chapter}, we will describe the RB method, the proper orthogonal decomposition (POD) \cite{pod_tol_error} used to construct the reduced space and the affine decomposition \cite{rozza_rom} exploited to ensure the efficiency of the method. Subsequently, in section \ref{cd_chapter}, the continuation method and the deflation one will be discussed with a particular focus on their implementation. The results obtained with the described methods will be presented in section \ref{results_chapter}, that will be divided as follows: in the first part we will focus on a bifurcation diagram with a single parameter to prove that it can be computed during the online phase whereas, in the second part, we will highlight the efficiency of the method considering an additional parameter. Finally, we will talk about some of the future perspectives of this work in section \ref{concl_chapter}.\\
We thus highlight that the advantages of the proposed method are strongly related to the way in which the four techniques are connected. In fact, even if the continuation and the deflation techniques allow one to compute a bifurcation diagram, they are very expensive and, without a reduced order model, the computational cost of the process may be prohibitive. Therefore, we decided to rely on the SEM and on the computation of a limited number of snapshots to perform the offline phase as efficiently as possible and to exploit the deflated continuation both in the offline and in the online phase. This way we could effectively obtain all the required solutions without exploiting any prior knowledge about the structure of the solution manifold.\\
We remark that the SEM is based on the open source software \textit{Nektar++} version 4.4.0 (see \cite{nektar}), while the reduced order model and the deflated continuation method described in this work have been implemented in \textit{ITHACA-SEM} (https://github.com/mathLab/ITHACA-SEM).

\section{Problem formulation}
\label{prob_form}
Let us consider the steady and incompressible Navier-Stokes equations \cite{ns} in the open bounded domain $\Omega\subset\R^2$ with a suitably regular boundary $\partial\Omega$:
\begin{equation} 
\label{NS_eq}
\begin{cases}
(\uu\cdot\nabla)\uu -\nu \Delta \uu + \nabla p=0 \quad &\text{in} \ \Omega, \\
\nabla \cdot \uu  =  0 \quad &\text{in} \ \Omega, 
\end{cases}
\end{equation}
where $\uu$ is the velocity, $p$ the pressure normalized over a constant density and $\nu$ the kinematic viscosity assumed constant. In order to better characterize the flow regime, it is important to observe that its main features can be summarized in a non-dimensional quantity named Reynolds number \cite{tritton} and defined as follows:
$$Re=\dfrac{UL}{\nu},$$
where $U$ is a characteristic velocity of the flow and $L$ is a characteristic length of the domain. System \eqref{NS_eq} can be obtained assuming that the Reynolds number is moderate and that the flow is steady.\\
It is important to observe that, when the equations in system \eqref{NS_eq} are normalized, they can be written in terms of non-dimensional quantities as:
\begin{equation} 
\label{NS_eq2}
\begin{cases}
(\uu\cdot\nabla)\uu - \frac{1}{Re} \Delta \uu + \nabla p=0 \quad &\text{in} \ \Omega, \\
\nabla \cdot \uu  =  0 \quad &\text{in} \ \Omega.
\end{cases}
\end{equation}
Since the structure of the mass balance equation does not change, while the only non-dimensional parameter in the momentum balance equation is the Reynolds number, one can conclude that it is the only meaningful one. Therefore, all the features present in the bifurcation diagrams that will be shown can be discussed in terms of $Re$. However, we decided to use as parameters the viscosity $\nu$ and a multiplicative factor on the  Dirichlet inlet boundary condition, which we will denote by $s$. It will be possible, in Section \ref{results_chapter}, to observe that the strict relation between the involved parameters and the Reynolds number allows one to explain the fact that the bifurcation points can be grouped, according to their nature, in specific and predictable curves.\\
To consider a well-posed problem, we supplemented system \eqref{NS_eq} with proper boundary conditions: a stress free boundary condition on the velocity at the outlet, a no-slip Dirichlet boundary condition on the physical walls and the following Dirichlet boundary conditions at the inlet
\begin{equation}
\uu  = \begin{bmatrix} u\\ v \end{bmatrix} = \begin{bmatrix} 20s(5-y)(y-2.5) \\ 0 \end{bmatrix} , \hspace{1cm} x=0,\text{ }y\in(2.5,5),
\label{real_bc}
\end{equation} 
where $s$ is the second parameter that we included in the model and a parabolic profile has been imposed to consider a more realistic condition. It should be observed that the dimension of the domain, described in figure \ref{mesh}, and the constant in front of the parabolic profile in the inlet boundary condition are used to obtain values of $Re$ approximately in $(60,260)$ 
when $\nu\in[0.3,1]$ and $s\in[0.8,1]$. We selected this range because, with the described geometry, it contains two bifurcations. Note that small changes in the geometry can influence the critical points positions, as observed in figure \ref{3d_diagr} or in \cite{martin_curve_walls}, it is thus important to find the bifurcation points associated with the exact geometric setting in order to apply the technique to real scenarios. \\
Furthermore, we introduce the variational formulation \cite{segal}, required by the Spectral element method (SEM) to obtain a discrete solution of the Navier-Stokes system. This will also be fundamental for the applicability of the Reduced basis method (RB). When deriving the weak formulation, one has to set appropriate functional spaces:
\[
\left(H_{0,\partial\Omega_D}^1(\Omega)\right)^2= \{\ww\in \left(H^1(\Omega)\right)^2 \mid \ww=0 \text{ on }\partial\Omega_D\},
\]
\[
L_0^2(\Omega)= \{q\in L^2(\Omega) \mid \int_\Omega qdx=0\},
\]
where $\partial\Omega_D$ is the portion of the boundary where Dirichlet boundary conditions are imposed.\\
Moreover, one multiplies the equations in system \eqref{NS_eq} by the appropriate test functions, integrates over the entire domain and integrates by parts the integrals associated with $\Delta\uu$ and $\nabla p$ obtaining the following weak formulation: find $\uu \in\uu^D+ (H_{0,\partial\Omega_D}^1(\Omega))^2$ and $p \in L_0^2(\Omega)$ such that
\begin{equation}
\label{gal_ns}
\left\{
\begin{aligned}
 \nu\int_\Omega\nabla \uu\cdot\nabla \vv +\int_\Omega \left((\uu\cdot\nabla)\uu\right)\cdot\vv - \int_\Omega p\nabla\cdot \vv  = 0 \quad &\forall \vv \in \left(H_{0,\partial\Omega_D}^1(\Omega)\right)^2, \\
\int_\Omega q\nabla\cdot \uu = 0\quad &\forall q \in L_0^2(\Omega),
\end{aligned}
\right.
\end{equation}
where $\uu^D\in H^1(\Omega)$ is a lifting function used to impose the non-homogeneous Dirichlet boundary conditions. Denoting $(H_{0,\partial\Omega_D}^1(\Omega))^2$ as $V$ and $L_0^2(\Omega)$ as $Q$, and introducing the following bilinear and trilinear forms:
\[
\begin{aligned}
a(\vv,\ww) =\nu\int_\Omega\nabla \vv\cdot\nabla \ww &\hspace{1cm}&&\forall\vv,\ww\in V,\\
b(\ww,q) = -\int_\Omega(\nabla\cdot \ww) \hspace{.05cm}q&\hspace{1cm}&&\forall\ww\in V,\text{ }\forall q\in Q,\\
c(\uu,\vv,\mathbf w)=\int_\Omega \left((\uu\cdot\nabla)\vv\right)\cdot\ww&\hspace{1cm}&& \forall \uu,\vv,\ww\in V,
\end{aligned}
\]
problem \eqref{gal_ns} can be expressed, in a more compact way, as follows:
find $\uu\in \uu^D+ V$ and $p\in Q$ such that
\begin{equation}
\label{gal_ns2}
\begin{cases}
 a(\uu,\vv) +c(\uu,\uu,\vv) +b(\vv,p)  = 0 \quad &\forall \vv \in V, \\
b(\uu,q) = 0\quad &\forall q \in Q.
\end{cases}
\end{equation}
Such a notation will be useful in section \ref{scm} to describe the static condensation method and we are going to exploit it to easily describe the two most common linearization techniques \cite{burger_notes}. See \cite{canuto_spettrali} for a more in-depth analysis of system \eqref{gal_ns2}.\\
The first one, named Oseen iteration, relies on the fact that, when the iterative solver is converging, two subsequent approximations are very similar. Therefore, the nonlinear term $c(\uu,\uu,\vv)$ can be approximated as follows:
\[
c(\uu^{k+1},\uu^{k+1},\vv)\approx c(\uu^{k},\uu^{k+1},\vv),
\]
where $\uu^k$ is the approximation obtained in the last iteration, while $\uu^{k+1}$ is the unknown one. Such a linearization technique is very common because it can be easily implemented and the associated iterative solver is very stable, but the convergence is only linear. On the other hand, one can exploit the Newton method to obtain a quadratic convergence, although a more accurate initial guess is required. To derive the latter, one expresses the unknown approximation as:
\[
\uu^{k+1}=\uu^k+\delta\uu,
\]
where $\delta\uu$ is the variation between the unknown solution $\uu^{k+1}$ and the last computed one $\uu^k$. This way, it is possible to approximate the nonlinear term $c(\uu,\uu,\vv)$ as follows:
\[
c(\uu^{k+1},\uu^{k+1},\vv)\approx c(\uu^{k+1},\uu^{k},\vv)+c(\uu^{k},\uu^{k+1},\vv)-c(\uu^{k},\uu^{k},\vv) . 
\]
In this work we exploited both techniques in order to increase the effectiveness of the deflation method (see Section \ref{deflation}). In fact, the Oseen iteration ensures a too slow convergence but, using only the Newton method, the iterative solver often diverges when solving the deflated problem.

\section{The spectral element method}
\label{sem_chapter}
In this section the main features of the SEM \cite{canuto_spettrali} will be described; since we are interested in the efficiency of the method, the main focus will be on the static condensation method. This is a technique that can be used to significantly reduce the computational cost to get the solution of the obtained linear system.\\
Let us consider the following Galerkin formulation \cite{canuto_spettrali2}, derived from problem \eqref{gal_ns2}: find $\uu\in\uu^D_{SEM}+ V^\delta$ and $p\in Q^\delta$ such that:
\begin{equation}
\label{gal_ns3}
\begin{cases}
 a(\uu,\vv) +c(\uu,\uu,\vv) +b(\vv,p)  = 0 \quad & \forall \vv\in V^\delta,\\
b(\uu,q) = 0 \quad & \forall q\in  Q^\delta,
\end{cases}
\end{equation}
{where $\uu^D_{SEM}$ is a suitable discretization of $\uu^D$ and $V^\delta$ and $Q^\delta$ are, respectively, two finite dimensional subspaces of $V$ and $Q$. It is thus possible to consider the bases $\{\phi_{u,i}\}_{i=1}^{N_u}$ and $\{\phi_{v,i}\}_{i=1}^{N_v}$ associated with the two components of the velocity and $\{\phi_{p,i}\}_{i=1}^{N_p}$ associated with the pressure. Therefore, discrete velocity and pressure can be expressed as follows:
\[
u=\sum_{i=1}^{N_u}u_i\phi_{u_i},
\hspace{1.5cm}
v=\sum_{i=1}^{N_v}v_i\phi_{v_i},
\hspace{1.5cm}
p=\sum_{i=1}^{N_p}p_i\phi_{p_i},
\]
where $u_i$, $v_i$ and $p_i$ are scalar coefficients that characterize the velocity and pressure fields.
Moreover, in the SEM the basis functions $\phi_{\cdot,i}$ are high-order polynomials inside the associated elements \cite{hp_fem}. In particular, in this work we decided to use the stable pair $\mathcal P_{P}(\Omega^e)/\mathcal P_{P-2}(\Omega^e)$, i.e. the velocity is represented by a polynomial of order $P$ while the pressure by one of order $P-2$ inside each element $\Omega^e$ (see \cite{p-p-21} and \cite{p-p-22} for a more detailed explanation of the approach). To compute the numerical results of section \ref{results_chapter}, even though this choice does not significantly influence them, we used $P=8$ when only 1 or 3 branches are involved and $P=12$ when the entire diagram is obtained. These values are chosen in order to balance the efficiency of the solver (the lower the order, the lower the computational cost) and its accuracy and the effectiveness of the deflation. In fact we observed that if the full order solver was more accurate, then the asymmetrical branches could be found more easily and earlier both in the offline and in the online phase.\\
The high order of the polynomials implies two main consequences: firstly, since several degrees of freedom are associated with each element, it is possible to obtain accurate solutions even with very coarse meshes like the one showed in figure \ref{mesh}. In general, comparing the SEM with the FEM, this approach is convenient because assembling and solving the linear system is expensive, but the computational cost of such operations is reduced if the solution is represented with fewer degrees of freedom, as in the SEM. Moreover, even if to compute bifurcation diagrams the mesh generation cost is negligible because one computes numerous solutions on the same mesh, generating fine meshes as the ones required by the FEM is an expensive operation \cite{mesh_refinement}. Secondly, this way it is possible to ensure the exponential convergence of the method instead of the algebraic one that characterizes the FEM \cite{hp_refinement} when the solution is smooth enough. This is important because the same level of accuracy can be obtained with a significantly lower number of degrees of freedom and, therefore, the assembly of the linear system and the computation of the solution are much more efficient. Furthermore, exploiting the high order polynomials, it is possible to rely on more accurate differentiation and integration formulas \cite{gauss_quadrature}. 

\subsection{The static condensation method}
\label{scm}
The static condensation method is a technique that allows one to solve a linear system much more efficiently, exploiting a specific structure of the associated matrix. It should be noted that, even if such a technique may be used to exploit a domain decomposition approach in a reduced basis framework \cite{sc_rb}, in this work we used it only in the offline phase to increase the efficiency of the SEM and, therefore, we did not link it with the RB method. Moreover, we remark that we will outline the method as described in \cite{nektar_online} (we are interested in the description of the method 
\textit{``void Nektar::CoupledLinearNS::SetUpCoupledMatrix''}).\\
To obtain the required structure, while solving the Navier-Stokes equations with the SEM, one has to split the velocity degrees of freedom into different groups. The first one contains the interior modes, i.e. all the basis functions with support inside a single element, while the second one contains all the remaining basis functions, that will be denoted as boundary modes \cite{sem}. It is crucial to observe that two interior modes associated with two different elements are orthogonal to each other because the measure of the intersection of their supports is zero. Such a property can be exploited to properly sort all the degrees of freedom to obtain block matrices with very small blocks associated with the elements. Denoting as $\uu_{int}$ the vector of the velocity degrees of freedom associated with the interior modes, as $\uu_{bnd}$ the vector associated with the element boundary ones and as $p$ the one associated with the pressure modes, it is possible to expand the linear system associated with problem \eqref{gal_ns3} as follows:
\begin{equation}
\label{nssc_system}
\begin{bmatrix}
A & -D_{bnd}^T & B \\
-D_{bnd} & 0 & -D_{int} \\
\widetilde B^T & -D_{int}^T & C \end{bmatrix}
\begin{bmatrix}
\uu_{bnd} \\ p \\ \uu_{int} \end{bmatrix} = \begin{bmatrix}\text f_{bnd} \\ 0 \\ \text f_{int}  \end{bmatrix}.
\end{equation}
It can be observed that the elements of the submatrix $A$ represent the relations between pairs of boundary modes, while the ones of the submatrices $B$ and $C$ are, respectively, associated with boundary-interior and interior-interior pairs. Finally, the interaction between the pressure and the velocity is summarized in $D_{bnd}$ and $D_{int}$. The elements of such matrices can be computed, $\forall i, j=1,\dots,N^{bnd}, \forall n, m=1,\dots,N^{int}$ and $\forall l=1,\dots,N^{p}$ as follows:
\begin{gather*}
A[i][j]=c(\boldsymbol\phi_{bnd}^j,\uu^k,\boldsymbol\phi_{bnd}^i)+c(\uu^k,\boldsymbol\phi_{bnd}^j,\boldsymbol\phi_{bnd}^i)+a(\boldsymbol\phi_{bnd}^j,\boldsymbol\phi_{bnd}^i), \\
B[i][n]=c(\boldsymbol\phi_{int}^n,\uu^k,\boldsymbol\phi_{bnd}^i)+c(\uu^k,\boldsymbol\phi_{int}^n,\boldsymbol\phi_{bnd}^i)+a(\boldsymbol\phi_{bnd}^n,\boldsymbol\phi_{int}^i), \\
\widetilde B^T[n][i]=c(\boldsymbol\phi_{bnd}^i,\uu^k,\boldsymbol\phi_{int}^n)+c(\uu^k,\boldsymbol\phi_{bnd}^i,\boldsymbol\phi_{int}^n)+a(\boldsymbol\phi_{bnd}^i,\boldsymbol\phi_{int}^n), \\
C[n][m]=c(\boldsymbol\phi_{int}^m,\uu^k,\boldsymbol\phi_{int}^n)+c(\uu^k,\boldsymbol\phi_{int}^m,\boldsymbol\phi_{int}^n)+a(\boldsymbol\phi_{int}^m,\boldsymbol\phi_{int}^n), \\
D_{bnd}[l][i]=b(\boldsymbol\phi_{bnd}^i,\phi_p^l), \\
D_{int}[l][n]=b(\boldsymbol\phi_{int}^n,\phi_p^l), \\
\text f_{bnd}[i]=f(\boldsymbol\phi_{bnd}^i)+c(\uu^k,\uu^k,\boldsymbol\phi_{bnd}^i), \\
\text f_{int}[n]=f(\boldsymbol\phi_{int}^n)+c(\uu^k,\uu^k,\boldsymbol\phi_{int}^n),
\end{gather*}
where $N^{bnd}$, $N^{int}$ and $N^p$ are, respectively, the number degrees of freedom associated to the velocity boundary modes, to the velocity interior ones and to the pressure ones. Moreover, the terms $f(\boldsymbol\phi_{bnd}^i)$ and $f(\boldsymbol\phi_{int}^i)$ are associated with the external forces that act on the system. However, since we assumed their absence, these two terms can be neglected.\\
It should be observed that such expressions hold when the Newton method is employed, however, when using the Oseen one, the first terms in the expansions of $A[i][j]$, $B[i][n]$, $\widetilde B^T[n][j]$ and $C[n][m]$ and the last ones in $\text f_{bnd}[i]$ and $\text f_{int}[n]$ have to be discarded.
Moreover, it should be noted that the submatrix $C$ is block diagonal and, therefore, it is easy to invert. Denoting $I$ as the identity matrix, we can premultiply the previous system by the matrix:
\begin{equation*}
\label{nssc_system2}
K = \begin{bmatrix} I & 0 & -BC^{-1} \\
 0 & I & D_{int}C^{-1} \\
 0 & 0 & I 
 \end{bmatrix},
\end{equation*}
obtaining the following one
\begin{equation}
\label{nssc_system3}
{\scriptstyle\begin{bmatrix} A-BC^{-1}\widetilde B^T & -D_{bnd}^T+BC^{-1}D_{int}^T & 0 \\
-D_{bnd}+D_{int}C^{-1}\widetilde B^T & -D_{int}C^{-1}D_{int}^T & 0\\
\widetilde B^T & -D_{int}^T &C \end{bmatrix}
\begin{bmatrix} \uu_{bnd} \\ p \\ \uu_{int}  \end{bmatrix} = \begin{bmatrix}\text f_{bnd}-BC^{-1}\text f_{int} \\ D_{int}C^{-1}\text f_{int} \\ \text f_{int}  \end{bmatrix}.}
\end{equation}
It is important to observe that, this way, the third equation has been decoupled from the other ones and the associated unknowns can be easily obtained after having solved the remaining $2\times 2$ block. Let us focus on the first $2\times 2$ block system involving $\uu_{bnd}$ and $p$, that can be written, simplifying the notation and considering a new set of matrices, as:
 \begin{equation}
\label{nssc_system4}
\begin{bmatrix}\hat A & \hat B \\ \hat C & \hat D \end{bmatrix}
\begin{bmatrix}\mathbf b \\ \hat p \end{bmatrix} = \begin{bmatrix} \hat{\text f}_{bnd} \\ \hat {\text f}_p  \end{bmatrix},
\end{equation}
where $\mathbf b=[\uu_{bnd},p_0]$ is a vector that contains both $\uu_{bnd}$ and the mean pressure $p_0$ (or a degree of freedom associated with it) while $\hat p$ accounts for the remaining pressure coefficients. A second level of static condensation can be obtained repeating the previous steps with the matrix:
\[
\hat K = \begin{bmatrix} I &\hspace{0.1cm} -\hat B \hat D^{-1} \\
 0 & I
 \end{bmatrix},
\]
in order to modify equation \eqref{nssc_system4} into the following one:
\begin{equation}
\label{nssc_system6}
\begin{bmatrix} \hat A-\hat B\hat D^{-1}\hat C & 0\\ \hat C & \hat D \end{bmatrix} 
\begin{bmatrix}\mathbf b \\ \hat p \end{bmatrix} = \begin{bmatrix}\hat {\text f}_{bnd}-\hat B\hat D^{-1}\hat {\text f}_p \\ \hat {\text f}_p  \end{bmatrix}.
\end{equation}
This way, one can solve the first row of equation \eqref{nssc_system6} to obtain $\mathbf b$, then one can substitute it into its second row to obtain $\hat p$ and, finally, these quantities can be used with the third row of problem \eqref{nssc_system3} to compute $\uu_{int}$. The computational cost of such a process is significantly lower than the direct computation of the solution of problem \eqref{nssc_system}, in fact $\mathbf b$ is the only unknown vector that has to be computed solving a linear system.

\section{The reduced basis method}
\label{rb_chapter}
In this section the RB method \cite{rozza_rom} will be briefly described. Such a technique can be used to efficiently compute the solution of a system of PDEs and, therefore, it is often used in optimization, real time queries, optimal control, design and uncertainty quantification. In this work, since we are interested in the discretization of bifurcation diagrams with many parameters, it has been used to compute the numerous required solutions. It should be observed that the same result could be obtained without exploiting the RB method, but the computational cost would be prohibitive in this case. Let us consider the following abstract parametric problem: given $\mu\in P$, find $u(\mu)\in \widetilde V$ such that
\begin{equation}
\label{fo_problem}
a(u(\mu),v;\mu)=f(v;\mu) \quad \forall v \in\widetilde V,
\end{equation}
where $\widetilde V$ is a suitable Hilbert space, $\mu \in P\subset \R^N$ is a vector of scalar parameters, $a(\cdot,\cdot;\mu):\widetilde V \times\widetilde V\rightarrow\R$ is a symmetric, coercive, bilinear and continuous operator for any parameter $\mu$ in the parameter space $P$. Analogously, $f(\cdot;\mu):\widetilde V\rightarrow\R$ is a linear and continuous operator for any $\mu \in P$. It is then possible to consider the discrete version of problem \eqref{fo_problem}, that reads as follows: given $\mu \in P$, find $u^\delta(\mu) \in\widetilde V^\delta$ such that
\begin{equation}
\label{sem_problem}
a(u^\delta(\mu),v;\mu)=f(v;\mu) \quad \forall  v \in\widetilde V^\delta,
\end{equation}
where $\widetilde V^\delta$ is a finite dimensional subspace of $\widetilde V$ that depends on the chosen discretization method. We remark that similar discrete and parametric problems have been deeply investigated, we thus refer to \cite{canpar1}, \cite{canpar2} and \cite{canpar3} for a more detailed analysis.\\
In this work the employed discretization method is the SEM, we can thus obtain accurate solutions with a space $\widetilde V^\delta$ smaller than the one used in the FEM. However, in order to significantly increase the efficiency of the solver, one would like to use spaces described by only few basis functions. Defining as $\widetilde V^{rb}\subset \widetilde V^\delta$ the finite-dimensional space used in the online phase of the reduced basis method, and assuming that $N^{rb}=\text{dim}(\widetilde V^{rb})\ll \text{dim}(\widetilde V^\delta)=N^\delta$, one can exploit the same formulation of problem \eqref{sem_problem} to define a reduced problem that can be solved much more efficiently: given $\mu\in P$, find $u^{rb}(\mu) \in\widetilde V^{rb}$:
\begin{equation}
\label{rb_problem}
a(u^\delta(\mu),v;\mu)=f(v;\mu) \quad \forall v \in\widetilde V^{rb}.
\end{equation}
It is thus necessary to be able to generate a small space $\widetilde V^{rb}$ capable to accurately discretize the continuous solutions of problem \eqref{fo_problem}. To do so, one can exploit problem \eqref{sem_problem} and different numerical techniques. In this work we only used the Proper Orthogonal Decomposition (POD) \cite{rozza_rom} to construct such a reduced space because its downsides are balanced by the advantages of the deflated continuation method \cite{tau_deflation} that will be described in Section \ref{cd_chapter}. In fact, two of the main disadvantages of the POD method are related to the high computational cost of the offline phase (indeed many solutions have to be computed) and to the fact that the parameter space has to be properly sampled. On the other hand, the deflated continuation method allows one to efficiently compute the required solutions and to automatically select only proper values of the parameter.\\
Finally, to link this section with the previous one, it is possible to consider $a(u^\delta(\mu),v;\mu)$ and $f(v;\mu)$ as the left hand side and the right hand side of both equations of problem \eqref{gal_ns2}, where $\uu^D$ is substituted by a suitable lifting function $u^{\delta,D}$ and $\widetilde V$ and $\widetilde V^\delta$ can be identified with $V\times Q$ and $V^\delta\times Q^\delta$ respectively. This way the operator on the left hand side is only a continuous tri-linear operator and the equation is stable if the inf-sup condition is satisfied. Even tough, it is convenient to enrich the reduced space with suitable functions named supremizers \cite{supremizer} to ensure such a stability condition, we observed that for moderate Reynolds numbers the described approach was already stable. Therefore, we did not enrich $\widetilde V^{rb}$ with the supremizers.

\subsection{The proper orthogonal decomposition}
\label{pod_section}
The POD can be used to generate a reduced space $\widetilde V^{rb}$ that is optimal in the euclidean norm over the space spanned by the snapshots. We are thus interested in the $\ell^2$ norm of the vectors containing the degrees of freedom with respect to the SEM basis, such a norm will be denoted as $\Vert\cdot\Vert$ in this section. \\
Let $P_M$ be a finite sampling of $P$ of dimension $M$, i.e. $P_M=\{\mu_1,...,\mu_M\}$. In order to construct $\widetilde V^{rb}$ one considers a symmetric and linear operator $C:\widetilde V^\MM\rightarrow\widetilde V^\MM$ defined as follows:
\[
C(v^\delta)=\frac 1M\sum_{i=1}^{M}\left(v^\delta,\psi_i\right)_{\widetilde V}\psi_i, \hspace{1.5cm} v^\delta\in\widetilde V^\MM,
\]
where $\psi_i=u^\delta(\mu_i)$ and $\widetilde V^\MM=\text{span}\{u^\delta(\mu):\mu\in P_M\}$. \\
Then, one computes the eigenvalue-eigenvector pairs $(\lambda_i,\xi_i)\in\R\times\widetilde V^\MM$, such that $\Vert \xi_i\Vert=1$ for any $i=1,...,M$, that satisfy $C\left(\xi_i\right)=\lambda_i\xi_i$. Finally, it is possible to sort the eigenvalues in descending order and the associated eigenvectors accordingly, and generate the reduced space with the first $N^{rb}$ eigenvectors. This way, it is possible to prove that the error obtained approximating the solutions of $\widetilde V^\MM$ with the ones in $\widetilde V^{rb}$ is given by:
\begin{equation}
\label{pod_error1}
{\frac{1}{M}\sum_{i=1}^{M}\Vert \psi_i-P_{N^{rb}}\left(\psi_i\right)\Vert^2}={\sum_{i=N^{rb}+1}^M\lambda_i},
\end{equation}
where $P_{N^{rb}}:\widetilde V^\MM\rightarrow\widetilde V^{rb}$ is a projection operator over $\widetilde V^{rb}$ defined as 
\[
P_{N^{rb}}(v)=\sum_{i=1}^{N^{rb}}\left(v,\xi_i\right)_{\widetilde V}\xi_i.
\]
Furthermore, it is important to observe that $V^{rb}$ is the only $N^{rb}$-dimensional space that minimizes the following quantity (see \cite{galerkin_POD} or \cite{pod_tol_error}):
\[
\mathop{\inf_{v^{rb}\in\widetilde V^{rb}\subset\widetilde V^\MM}}_{\text{dim}(\widetilde V^{rb})=N^{rb}} 
 \sum_{i=1}^M \Vert \psi_i-v^{rb}\Vert^2,
\] 
and that, to obtain a numerically stable online solver, it is convenient to orthonormalize the eigenvectors.\\
Therefore, from the implementation point of view, in order to generate a reduced space with the desired properties, one has to compute the eigenvalues and the orthonormalized eigenvectors of the correlation matrix of the snapshots and has to express the latter as linear combinations of full order basis functions. One of the most efficient method to perform such a task is the Singular Value Decomposition (SVD) \cite{fast_svd}. It is additionally important to note that, when interested in generating a space that is optimal with respect to a different norm, one has to premultiply the involved correlation matrix with the Cholesky factor of the matrix associated with the corresponding inner product before performing the described operations \cite{weighted_POD}.

\subsection{The affine decomposition}
\label{ad_section}
In order to efficiently solve problem \eqref{rb_problem}, it is important to seek the solution in a lower dimensional space, as the one that can be obtained through the POD method, whose accuracy is ensured thanks to property \eqref{pod_error1}. However, in order to benefit of this low dimensionality, by splitting the computation in the two phases, the model has to fulfill some additional hypothesis. In particular, in order to have an online phase which is independent of the number of degrees of freedom of the approximation, namely $N^\delta$, one can exploit the so-called affine decomposition \cite{supremizer2}.\\
The main idea is to precompute several matrices and vectors in the offline phase that will be used for every instance of a new parameter, during the online one, to rapidly assemble the linear system. This is important because, this way, such an assembly does not scale with the dimension of $\widetilde V^\delta$ and all the operations of the online phase only depend on $N^{rb}$. Let us denote the linear system that has to be solved in the online phase as
\[
A^{rb}(\mu)u^{rb}=\text f^{\hspace{0.05cm}rb}(\mu).
\]
Using this notation, the affine decomposition assumption reads as:
\[
A^{rb}(\mu)=\sum_{q=1}^{Q_a}\theta_a^q(\mu)A_q^{rb},
\]
\[
{\text f}^{\hspace{0.05cm}rb}(\mu)=\sum_{q=1}^{Q_f}\theta_f^q(\mu){\text f}_q^{\hspace{0.05cm}rb}.
\]
It is important to observe that $A^{rb}(\mu)$ and ${\text f}^{\hspace{0.05cm}rb}(\mu)$ are expressed as linear combinations of matrices or vectors that do not depend on the parameter and, therefore, its contribution is entirely included in the coefficients. This way it is possible to precompute $A^{rb}_q$ for any $q=1, \dots, Q_a$ and ${\text f}_q^{\hspace{0.05cm}rb}$ for any $q=1, \dots, Q_f$ in the offline phase and, subsequently, assemble $A^{rb}(\mu)$ and ${\text f}^{\hspace{0.05cm}rb}(\mu)$ computing only the scalar coefficients. However, it should be noted that it is not always possible to directly exploit such a structure but, when required, it can be approximated with the so called Empirical Interpolation Method (EIM) \cite{eim}. We remark that we did not use the EIM because we projected the linearized operators onto the reduced space to directly exploit the affine decomposition.

\section{Numerical computation of bifurcation diagrams}
\label{cd_chapter}
Let us consider the following nonlinear parametric equation:
\begin{equation}
\label{general_eq}
L(\uu;\mu)=0 ,
\end{equation}
where the function $\uu$ belongs to a suitable functional space $\overline V$, the parameter $\mu$ belongs to $\R^N$ and $L$ is a nonlinear operator. Note that, in this paper, we can consider $\overline V=V\times Q$ because we are interested in computing a diagram associated with problem \eqref{gal_ns2}.\\
Since, due to the nonlinearity, several solutions can exist for the same value of the parameter, it is important to summarize them in a single diagram in order to highlight the main properties of the solutions manifold associated with the parameter space. Such diagrams are called bifurcation diagrams \cite{kielhofer}, the information is often summarized by means of a scalar output function, while the parameter is represented on the other axis. It should be noted that, if $n>1$, more than one axis are required to properly represent the parameter and it could be useful to rely on more advanced visualization techniques to show the diagram \cite{bif_diagr_nd}.\\
In this work we decided to compute the bifurcation diagrams with a combination of the continuation method and the deflation method. The main idea behind the coupling of this two techniques is that, on the one hand, one can entirely discretize a specific branch of the diagram given a solution belonging to it while, on the other hand, the deflation method is used to compute the first solutions of the new branches that are required by the continuation method. Such an approach, where the two described techniques are alternated in order to discover and follow each branch of the diagram, is called deflated continuation and is described in \cite{tau_deflation}. However, we implemented a more elaborated version of it, in fact we decided to use two different versions of the continuation method and to pair a novel approach to the deflation one in order to increase the efficiency and the effectiveness of the deflated continuation method. Such modifications will be better explained in the following sections.

\subsection{Continuation method}
\label{continuation}
In order to simplify the notation in the discussion, let us assume that, in problem \eqref{general_eq}, there exists a solution $\uu$ for any value of the parameter $\mu$ and that $\mu\in\R$. The assumption $\mu\in\R$ is useful to simplify the description of the method and it is not restrictive, indeed one can consider $\mu\in\R^N$ and let vary only one component at a time. This way the discussed approach can be extended to problems characterized by more parameters.\\
Since problem \eqref{general_eq} is nonlinear, an iterative solver is required to compute a solution. However, in order to ensure its convergence, it needs an initial guess close enough to the sought solution. Such a guess can be obtained with the continuation method. The aim of the technique, in fact, is to compute a proper initial guess in order to allow the solver to converge to a solution on the same branch of the last one. This way, any arbitrary branch can be entirely reconstructed repeating several times the following procedure \cite{seydel}.\\
Initially, one assumes to know $m$ solutions ($\uu_1, \dots, \uu_m$) on a branch of the diagram, they are associated with the parameter values $\mu_1, \dots, \mu_m$ and one wants to compute the solution $\uu_{m+1}$ associated with the value $\mu_{m+1}$. Given the input $\{(\uu_1,\mu_1), \dots, (\uu_m,\mu_m)\}$, the output of the continuation method is a function $\widetilde \uu_{m+1}$ close to the unknown solution $\uu_{m+1}$. Using $\widetilde\uu_{m+1}$ as an initial guess, the iterative solver can efficiently converge to $\uu_{m+1}$. Then, the process can be repeated to compute $\uu_{m+2}$ exploiting $\uu_{m+1}$ and the previous solutions, and so on.\\
In this work we implemented two different versions of the continuation method. The first one, that we will denote as simple continuation, is the simplest way to exploit the already computed information to obtain an initial guess and is characterized by $m=1$. On the other hand, the pseudo-arclength continuation \cite{keller} is a more advanced technique and exploits the last two solutions to compute the subsequent initial guess.\\
Let us first consider the simple continuation. Exploiting the continuity of each branch, one can assume that, if the step size $\Delta\mu_{i+1}=\mu_{i+1}-\mu_{i}$ is very small, the solutions $\uu_{i}$ and $\uu_{i+1}$ will be very similar to each other. This way, $\uu_{i}$ can be considered a good approximation of $\uu_{i+1}$ and used as initial guess to compute it. The advantages of such an approach are that it is inexpensive and requires a single solution. However, properly setting the quantity $\Delta\mu_{i+1}$ is complex, even if it is possible to rely on prior knowledge or heuristics. \\
Unfortunately, such a choice is very important, in fact a too wide step may make the solver diverge or converge to solutions on other branches, because of the significant difference between $\uu_{i}$ and $\uu_{i+1}$. On the other hand, a too short step would ensure the convergence to the correct solution but the computational cost would dramatically increase \cite{allgower}.\\
In order to properly set the step size $\Delta\mu_{i+1}$ and to improve the accuracy of the obtained initial guess, one can choose to rely on the pseudo-arclength continuation method. In such a technique the next value of the parameter $\mu_{i+1}$ is considered as an unknown and an alternative parametrization of the branch, characterized by the curve arclength $S$, is taken into account. To derive the system that has to be solved, as explained in \cite{keller}, let us consider the following normalization equation:
\begin{equation}
\label{normalization_L_NS}
N(\uu,\mu;\Delta S_i) \doteq \dot \uu_i^T(\uu-\uu_i)+\dot \mu_i(\mu-\mu_i)-\Delta S_i = 0,
\end{equation}
where $(\uu_i,\mu_i)$ is a point on a regular portion of the branch $\mathcal C$ and $(\dot\uu_i,\dot\mu_i)$ is the unit tangent to the curve in such a point. Equation \eqref{normalization_L_NS} characterizes the plane orthogonal to the vector $(\dot \uu_i,\dot\mu_i)$ such that the distance between $(\uu_i,\mu_i)$ and its projection on the plane is $\Delta S_i$. Moreover, if the line described by $(\dot\uu_i,\dot\mu_i)$ is a good approximation of $\mathcal C$, the orthogonal projection of $(\uu_i,\mu_i)$ on the plane is very similar to the sought solution $(\uu_{i+1},\mu_{i+1})$. A representation of such a process is outlined in figure \ref{continuation_picture}. 
\begin{figure}
\centering
\includegraphics[height=7.5cm, width=7.5cm, angle=0, keepaspectratio]{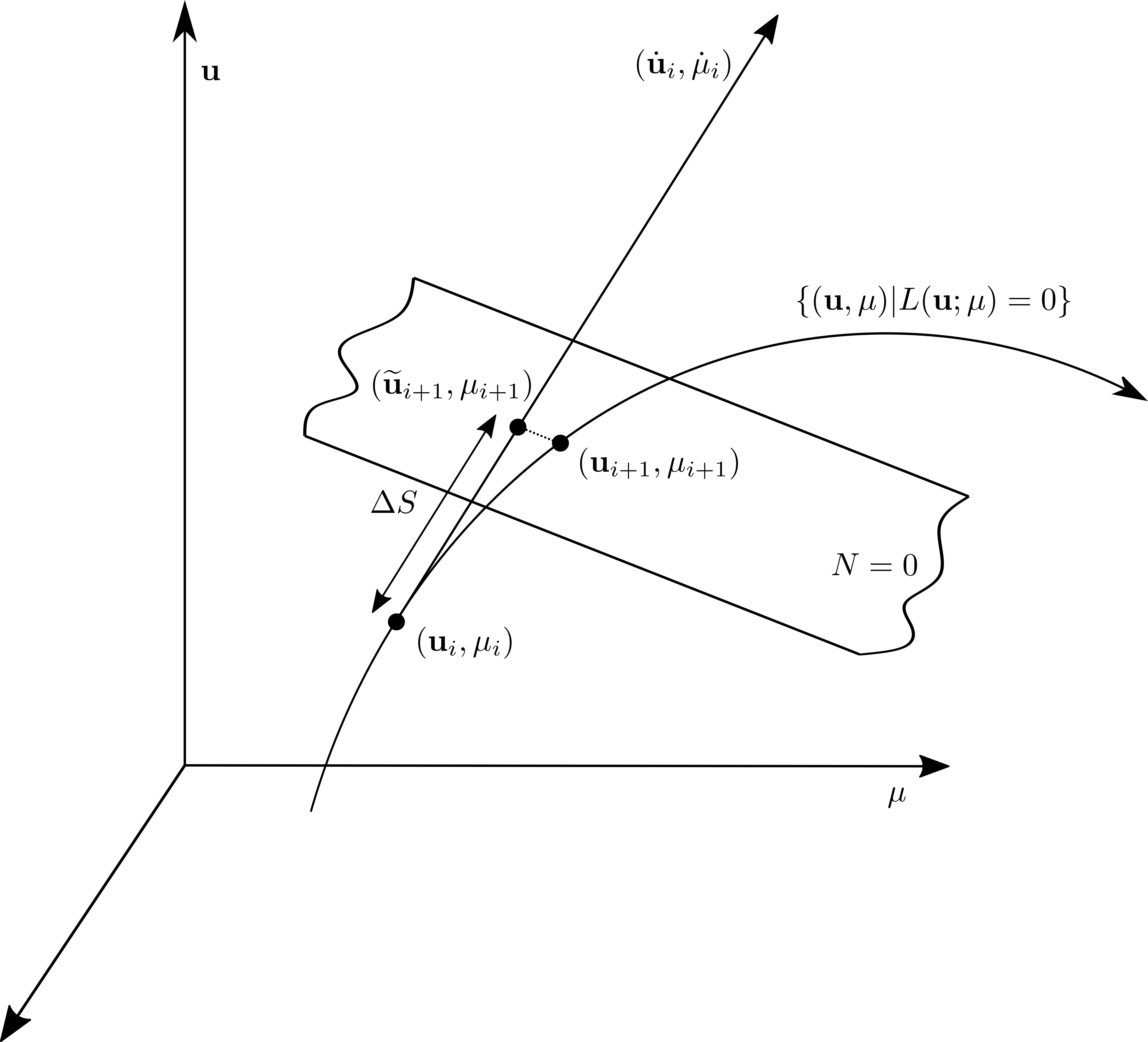}
\caption{Visualization of the pseudo-arclength continuation method}
\label{continuation_picture}
\end{figure}
Consequently, this projection can be used as a good initial guess. In order to compute it, one can solve the linear system:
 \begin{equation}
 \label{linear_continuation}
\begin{bmatrix} 
L_\uu^i & L_\mu^i \\
\dot \uu_i & \dot \mu_i 
\end{bmatrix}
\begin{bmatrix} 
\Delta \uu_i \\ \Delta \mu_i 
\end{bmatrix} = 
\begin{bmatrix} 
L^i \\ \Delta S_i 
\end{bmatrix} ,
\end{equation}
where the subscripts are associated with the derivation operation, the superscripts represent the point where the function is evaluated, i.e. $L_\uu^i=L_\uu(\uu_i;\mu_i)$ and $L_\mu^i=L_\mu(\uu_i;\mu_i)$, and the following notation has been used:
\[
\Delta \uu_i =\widetilde \uu_{i+1}-\uu_i,\hspace{2cm} \Delta \mu_i = \mu_{i+1}-\mu_i.
\]
This system can be obtained by linearizing the following one, obtained combining equation \eqref{general_eq} with equation \eqref{normalization_L_NS}, with the Newton method:
\begin{equation} \label{continuation_system}
\left\{
  \begin{aligned}
 &L(\uu;\mu)=0,\\
& \dot \uu_i^T(\uu-\uu_i)+\dot \mu_i(\mu-\mu_i)-\Delta S_i = 0 .
  \end{aligned}
  \right.
\end{equation}
Furthermore, since the quantities $\dot \uu$ and $\dot\mu$ are, in general, not available, they have to be approximated. We decided to use the following approximations:
\[
\dot\uu\simeq\dfrac{\uu_i-\uu_{i-1}}{\Delta S_{i-1}}, \hspace{2cm} \dot\mu\simeq\dfrac{\mu_i-\mu_{i-1}}{\Delta S_{i-1}},
\]
even if several alternatives exist.\\
The main advantage of such a version of the continuation method is that the subsequent value of the parameter is automatically chosen. This way, the steps are wider in very smooth regions, while they can be much shorter near the singularities. This is important because, when one wants to compute a bifurcation diagram, there are regions where the solution varies very rapidly, and regions where two solutions are very similar even if they are associated with two values of the parameters very far from each other.\\
We also decided to iteratively modify $\Delta S_i$ to further improve the effectiveness of the method, even though good results can be obtained also fixing it after some experimental observations. We chose to increase $\Delta S_i$ each time the solver converged in less than 6 iterations and to decrease it otherwise as suggested in \cite{seydel} for a similar setting. In fact small values of $\Delta S_i$ imply that the plane described in equation \eqref{normalization_L_NS} is very close to $\uu_i$ and that the solver can easily converge to $\uu_{i+1}$ because it is very similar to $\uu_i$. However, with such a choice, the computational cost of the diagram computation increases because of the excessive number of computed solutions. On the other hand, if $\Delta S_i$ is too large, $\widetilde \uu_{i+1}$ is a poor approximation of $\uu_{i+1}$ because the plane is very too far from $\uu_i$, therefore the solver may converge in several iterations or may not converge.\\
In general, it is possible to observe that such regions are, respectively, close to and far from a bifurcation point. Moreover, the pseudo-arclength continuation is more accurate than the simple one because it relies on a branch linearization. However, sometimes such a technique cannot be used. This issue can arise in two different scenarios. Firstly, when one wants to compute the second solution, only a single solution is available and, therefore, the technique cannot be applied. Secondly, right after a bifurcation point, the number of solutions varies from one iteration to the next one. This implies that it is not possible to exploit two solutions on the same branch to compute the initial guess. In such scenarios, we decided to use the simple continuation with a step size proportional to the last value of $\Delta\mu_i$ after a bifurcation point and with a very short step after the computation of the first solution. \\
Finally, we highlight that, since the structure of the matrix associated with problem \eqref{linear_continuation} is different from the one of problem \eqref{nssc_system}, we decided to implement a bordering algorithm \cite{keller} to solve the former with the static condensation method.

\subsection{Deflation method}
\label{deflation}
As discussed in the previous section, the continuation method allows one to follow each branch of the diagram, however, it requires a first solution on the branch to properly work. In this work such solutions have been computed in two different ways. We decided to compute the very first solution using the zero initial guess because of the lack of prior knowledge, while we used the deflation method to obtain the first solutions on unknown branches. Such a technique has been initially developed to compute multiple roots of a polynomial and, in order to explain it, it is convenient to discuss, as in \cite{classic_deflation}, a simplified scenario first.
Let us consider the scalar polynomial
\begin{equation}
\label{polynomial}
p(x) =c_0 \prod_{j=0}^{m-1} (x-x_j),
\end{equation}
where $c_0\in\R$ is a scaling factor and $x_0, \dots, x_{m-1}$ are $m$ distinct roots. Moreover, let us assume to be able to numerically compute just a single solution for any arbitrary polynomial with a numerical method. This way, one can easily compute root $x_i$ but one cannot obtain the other ones. However, to overcome such an issue, it is possible to consider the deflated polynomial
\begin{equation}
\label{polynomial2}
p_1(x) =c_0 \dfrac{1}{(x-x_i)}\prod_{j=0}^{m-1} (x-x_j).
\end{equation}
It should be observed that the polynomial in equation \eqref{polynomial2} is characterized by the same roots of the one in \eqref{polynomial} except for the $i$-th one. It is therefore possible to exploit the available algorithm to obtain a root of $p_1(x)$, consequently obtaining the second one of $p(x)$. Such a process can be repeated in order to obtain all the existing roots of $p(x)$. However, it should be remembered that, when such roots are computed with a numerical method, the obtained results are approximations of the exact ones and, therefore, the deflated polynomial is still characterized by all the roots of the previous one. Anyway, if the approximation is accurate enough, the numerically deflated polynomial is very similar to the analytical deflated one and the difference is significant only in a very small neighbourhood of the deflated root. It is thus possible to exploit such a technique to compute distinct roots that are far from each other.\\
The described approach can be then generalized to be applied to systems of PDEs, however, before discussing the generalized method, it is convenient to introduce the concept of deflation operator. In equation \eqref{polynomial2} the function $1/(x-x_i)$ is called deflation operator and it is responsible for removing a root from the polynomial; such an operator can be generalized in the following way (see \cite{classic_deflation}).
\begin{definition}
\textit{Let us denote by $W$ and $Z$ two Banach spaces and by $U$ an open subset of the additional Banach space $\hat V$. Moreover, let $L:U\subset\hat V\rightarrow W$ be a Fr\'echet differentiable operator and $L'$ be its Fr\'echet derivative. Then, let $M(\uu;\ww):W\rightarrow Z$ be an invertible linear operator for each $\ww \in U$ and for each $\uu\in U\backslash\{\ww\}$. If, for any Fr\'echet differentiable operator $L$ for which the following properties hold:
\[
L(\ww)=0, \hspace{2cm}L'(\ww)\text{ is nonsingular},
\]
and for any arbitrary sequence $\{\uu_i\}\subset U\backslash \{\ww\}$ such that $\uu_i\rightarrow\ww$, the following inequality holds
\begin{equation}
\label{deflation_definition_formula}
\liminf_{i\rightarrow\infty}\Vert M(\uu_i;\ww)L(\uu_i)\Vert_Z>0,
\end{equation}
then $M$ is a deflation operator.}
\end{definition}
In order to generalize equation \eqref{polynomial2}, one considers the deflation operator $M(\uu;\ww)$ and the following deflated system:
\begin{equation}
\label{general_eq2}
G(\uu;\mu) \doteq M(\uu;\ww)L(\uu;\mu)=0,
\end{equation}
that is characterized by the same solutions of $L(\uu;\mu)=0$ except for $\ww$.\\
Denoting by $I$ the identity operator, it is possible to consider the most straightforward generalization of the deflation operator introduced in equation \eqref{polynomial2}:
\begin{equation*}
\label{def1}
M(\uu;\ww)=\dfrac{I}{\Vert \uu-\ww\Vert_U}.
\end{equation*}
It can be observed that, since such an operator tends to zero when $\uu$ is very far from $\ww$, the iterative solver may converge to unphysical solutions if it is directly employed because the exact residual is multiplied by a factor that tends to zero. Even if several alternatives exist, we preferred to implement a very simple deflation operator, that can be expressed as:
\begin{equation*}
\label{def2}
M(\uu;\ww)=I+\dfrac{I}{\Vert \uu-\ww\Vert_U^p},
\end{equation*}
with $p=1$ in the offline phase and $p=2$ in the online; such quantities have been fixed via experimental observations. Note that, when $p$ increases, the region of attraction of $\ww$ is wider and, therefore, it is not possible to converge to fields very close to it. We observed that, if a solver is less stable (as the used online solver without the supremizer stabilization), it is convenient to use a higher value of $p$ in order to avoid computing solutions on the same branch more than once.\\
The main advantage of the deflation method consists in the ability to discover unknown branches without any prior knowledge, however, if other branches exist, one cannot be sure that they will be found with such a technique. In fact, if a branch $\mathcal C$ is too far from any known solution, the solver may diverge before reaching the region of attraction of any solution in $\mathcal C$. Therefore, it is advisable to fix a meaningful maximum number of iterations for the iterative solver when the deflated system is solved. Such a threshold should be high enough in order to let enough time to the iterative solver, however, if too many iterations are available, a lot of computational resources will be wasted when new branches cannot be found. \\
We decided to fix such a quantity equal to 150 in the offline phase and to 300 in the online one (because each iteration is significantly less expensive). Consequently, the deflation method is the bottleneck of the deflated continuation when applied at each step because the linear system has to be assembled and solved many times. Therefore, if one knows the position of the bifurcation points, it is advisable to use the deflation only in those regions.\\
Moreover, as described in \cite{tau_deflation}, it is possible to increase the efficiency of the deflation observing that one does not need to assemble and solve the linear systems associated with the deflated problem. In order to describe a more efficient approach, let us assume that one wants to solve system \eqref{general_eq2} with the Newton method, that is characterized by the following iteration:
\[
J_G\Delta\uu_G=-G ,
\]
where $J_G$ is the associated Jacobian. In the same way, we will refer to the corresponding iteration on the undeflated system as:
\begin{equation}
\label{l_system}
J_L\Delta\uu_L=-L .
\end{equation}
Assuming that $M(\uu;\ww)\in\R$ and exploiting the Sherman-Morrison formula \cite{golub}, one can derive the following relation:
\begin{equation*} 
  \begin{aligned}
 \Delta\uu_G &=-J_G^{-1}G\\
 		&=-\left(M J_L+L{M'}^T\right)^{-1}(M L)\\
 		&=-\left(M^{-1}J_L^{-1}-\dfrac{M^{-1}J_L^{-1}L{M'}^TM^{-1}J_L^{-1}}{1+{M'}^TM^{-1}J_L^{-1}L}\right)\left(M L\right)\\
 		&=\left(1-\dfrac{M^{-1}{M'}^T J_L^{-1}L}{1+M^{-1}{M'}^T J_L^{-1}L}\right)\left(-J_L^{-1}L\right).
  \end{aligned}
\end{equation*}
One can observe that the quantity that multiplies the term $\left(-J_L^{-1}L\right)$ is a scalar quantity, therefore we will refer to it as
\[
\tau \doteq 1-\dfrac{M^{-1}{M'}^T J_L^{-1}L}{1+M^{-1}{M'}^T J_L^{-1}L}=\dfrac{1}{1+M^{-1}{M'}^T J_L^{-1}L} .
\]
Finally, if one is able to efficiently compute a solution of the undeflated system, it is possible to exploit the same algorithm to obtain a solution of the deflated one. In fact, one simply has to multiply each Newton iteration for $\tau$, that can be computed, after having solved system \eqref{l_system}, very efficiently with a scalar product and scalar operations. Its cost is thus linear with the number of degrees of freedom. On the other hand, if one wants to directly solve the deflated system, one has to explicitly construct the matrix $J_G$, that is full, and solve the associated linear system without exploiting the sparsity of the original matrix $J_L$.

\subsubsection{An approach to improve the deflation method}
\label{heuristic_section}
In the previous section we explained how to efficiently deflate a system modifying the residual of the iterative solver associated with the undeflated problem instead of constructing and solving the deflated one. However, we analyzed the values assumed by the new scalar factor $\tau$ and we observed that they were always very low (very often its absolute value was lower than $\left.10^{-7}\right)$. This phenomenon implied serious consequences because, using the formula
\[
\Delta\uu_G=\tau\Delta\uu_L,
\]
one observes that $\Vert\Delta\uu_G\Vert\approx0$ when $\tau\approx0$. Consequently, the iterative solver would have required too many iterations and the computational cost would have been prohibitive. We thus analyzed the behaviour of $\tau$ and its relation with the solutions exploiting the fact that, since $\tau\in\R$, it cannot change the direction of the undeflated residual. However, even if the direction cannot differ between $\Delta\uu_G$ and $\Delta\uu_L$, their orientation can because $\tau$ can be positive or negative. Therefore, we decided to modify the values of $\tau$ while maintaining its sign to avoid losing changes in the orientation.\\
Firstly, we decided to fix two lower bounds in order to prevent $|\tau|$ to assume values too close to zero through the following formula:
\begin{equation} \label{small_taus}
\tau=\left\{
  \begin{aligned}
 &\tau_t^-  &\text{if }\tau>\tau_t^- ,\\
 &\tau_t^+ &\text{if }\tau<\tau_t^+ ,\\
 &\tau &\text{otherwise .} \\
  \end{aligned}
  \right.
\end{equation}
In this work we decided to use $\tau_t^-=-0.4$ and $\tau_t^+=0.6$, however, such values are very problem dependent. 
Secondly, we multiplied $\tau$ by a scaling factor $c$ in order to avoid using the thresholds too often. It is important to observe that we initially fixed $c=1$ and then we multiplied it with a suitable scaling factor $s_c$ each time $\tau$ switched from negative to positive. In fact, this change of sign implies that the deflation was preventing the solver to converge to a solution $\uu_0$ ($\tau<0$) but, when the current iteration was too far from it, the solver was attracted again by $\uu_0$ ($\tau>0$).  Instead, we observed that bigger values of $c$ could help the solver to better escape from the region of attraction of a solution and, therefore, we slightly increased it each time such an escape failed.
\begin{equation*} \label{small_taus_c}
\tau=\left\{
  \begin{aligned}
 &\tau_t^- &\text{if }c\tau>\tau_t^- ,\\
 &\tau_t^+ &\text{if }c\tau<\tau_t^+ ,\\
 &c\tau &\text{otherwise .} \\
  \end{aligned}
  \right.
\end{equation*}
Finally, we decided to use the Newton iteration when the current iteration was close to known solutions, while we exploited the Oseen one otherwise to improve the stability of the solver.\\
In this work we used $\tau_t^-=-0.4$, $\tau_t^+=0.6$ and $s_c=1.75$. These values have been chosen in order to obtain the best possible convergence velocity for the problem of interest, nevertheless different choices would lead to the same results at the cost of a slower convergence. Therefore, if prior knowledge is available, it is convenient to use it to optimally set these three quantities, otherwise one could use a relatively small value of $\tau_t^-$ and $\tau_t^+$ (for instance $-0.1$ and $0.1$) and a scaling factor $s_c$ slightly greater than 1 (for example $1.1$). This non-optimal choice may lead to a less effective deflation because it would be harder to escape from the region of attraction of the known solution, but the new branches could still be obtained simply increasing the number of available iterations. On the other hand, if the used values are significantly greater than necessary, the position of the bifurcation points would be computed less accurately but the new branches could still be continued to generate the entire diagram. We thus suggest to use very small values if the goal is to accurately localise the bifurcation points, or relatively large values when one wants to find all the branches.

\section{Numerical results}
\label{results_chapter}
In this section we will discuss the results that can be obtained using a combination of the described techniques. Let us remember that, in order to obtain a bifurcation diagram during the online phase, we performed the following steps. Firstly, the SEM (see section \ref{sem_chapter}) has been used to compute the snapshots involved in the reduced space generation. In order to choose the proper values of the parameter $\mu$ and obtain multiple solutions for any $\mu$, we respectively exploited the continuation and the deflation methods (see section \ref{cd_chapter}). Subsequently, all the snapshots are grouped together to construct a global reduced space with the POD (see section \ref{rb_chapter}) and, finally, the continuation and the deflation are used again in the online phase to efficiently reconstruct the diagram. Such steps are summarized in algorithm \ref{pseudo_code}.
\begin{minipage}{1\linewidth}\begin{algorithm}[H]
\label{pseudo_code}
  \caption{Main steps to efficiently compute a bifurcation diagram. $S_i$ is, for any $i=2,...,N$, a small set of values used to compute the snapshots, while $P_{-\mu_1}$ represents the parameter space for the parameters $\mu_2,...,\mu_N$.}
  \Offline{}{
  snapshots\_set = $\emptyset$;\\
   \For{$(\mu_2,...,\mu_N)\in S_2\times...\times S_N$}{
    // Compute a one-dimensional diagram with the function \\// ``offline\_deflated\_continuation()''\\
    new\_snapshots = offline\_deflated\_continuation();\\
   snapshots\_set = snapshots\_set $\cup$ new\_snapshots;    
   }
   // Generate the reduced space with the POD\\
   Reduced\_space = POD(snapshots\_set);
  }
  \Online{}{
    solutions\_set = $\emptyset$;\\
   \For{$(\mu_2,...,\mu_N)\in P_{-\mu_1}$}{
    // Compute a one-dimensional diagram with the function \\//``online\_deflated\_continuation()''\\
    new\_solutions = online\_deflated\_continuation();\\
   solutions\_set = solutions\_set $\cup$ new\_solutions;    
   }  }\end{algorithm}\end{minipage}\\[0.25cm]
   Since one of the hypothesis of the reduced basis method is the presence of a smooth solution manifold with just a single solution associated with each parameter value that, in this work, does not hold because of the pitchfork bifurcation points, we will first prove that such an approach is able to accurately discretize a bifurcation diagram with a single parameter. For instance, in figure \ref{curved_man}, a bifurcation diagram where such hypothesis do not hold is shown. Subsequently, we will move on to the description of a bifurcation diagram with two parameters, where we will be able to appreciate the efficiency of the computation. We also remark that all the shown errors are computed in the $L^2$ norm.\\
\begin{figure}
\centering
\includegraphics[height=5.5cm, width=5.5cm, angle=0, keepaspectratio]{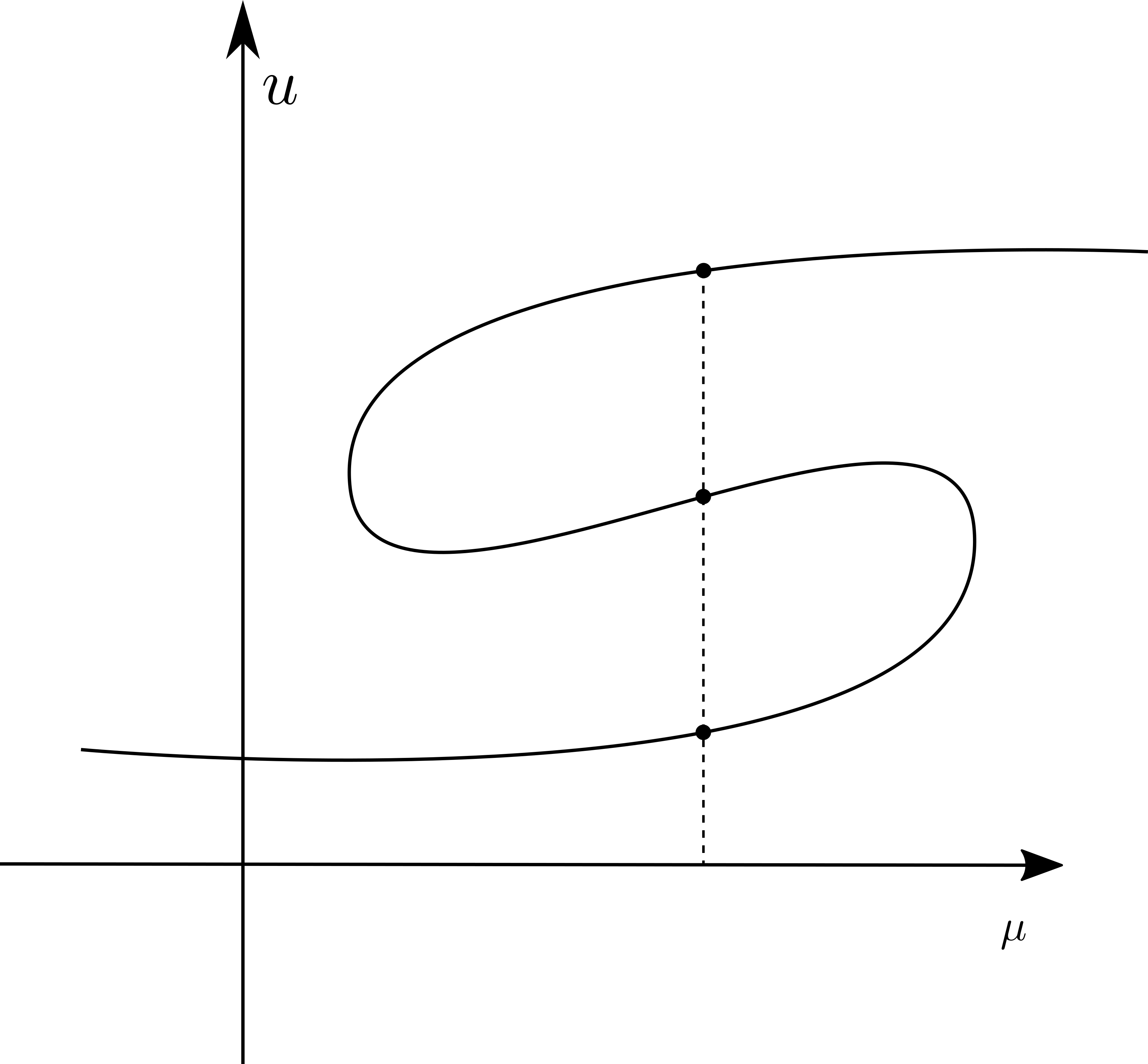}
\caption{Bifurcation diagram with multiple solutions associated with the same parameter value}
\label{curved_man}
\end{figure}In order to better understand the different existing solutions, it is convenient to analyze the ones in figure \ref{possible_solutions}, where the horizontal velocity is highlighted by the colour gradient. It is important to observe that more than one solution is associated with the values 0.6 and 0.3 of the viscosity. Moreover, some of them are axisymmetric while others are not. However, due to the symmetry of the domain and of the boundary conditions, it is always possible to reflect a solution over the horizontal axis of symmetry to obtain another solution. Such a phenomenon can be observed, for instance, in figures \protect{\subref{sol2}} and \protect{\subref{sol3}}. This is important because the function that will be used to compute the bifurcation diagram is, as suggested in \cite{classic_deflation}, the following one:
\[
f(\uu) \doteq \mypm\int_\Omega\Vert\uu-\mathcal R(\uu)\Vert^2,
\]
where $\uu$ is the velocity, $\mathcal R(\cdot)$ is an operator that reflects a solution over the horizontal symmetry axis and the sign is positive one if the jet hugs the upper wall or negative one otherwise. Therefore, $f(\uu)$ will be equal to zero if a solution is perfectly symmetrical, while its absolute value will increase with the asymmetry of the velocity field. Exploiting such an interpretation, one can immediately conclude that the bifurcation diagrams will be symmetrical because each solution $\uu_0$ can be mirrored over the symmetry axis obtaining another solution $\uu_1$ such that $f(\uu_0)=-f(\uu_1)$. Note that in the shown diagram the norm used to compute $f(\uu)$ is the $L^2$ norm, even though any norm could be used because we are only interested in obtaining different values of $f(\uu)$ for solutions on different branches.
\begin{figure}
\centering  
  \begin{subfigure}{5cm}
    \centering\includegraphics[width=5cm]{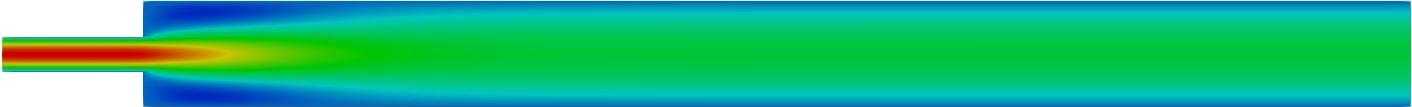}
    \caption{$\nu = 1$}
    \label{sol11}
  \end{subfigure}
  \begin{subfigure}{5cm}
    \centering\includegraphics[width=5cm]{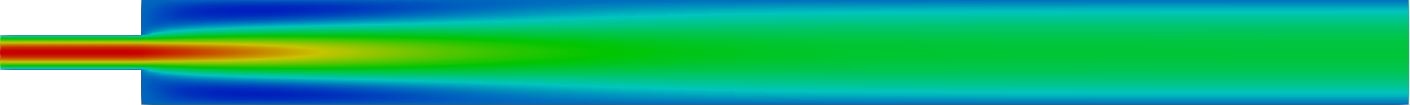}
    \caption{$\nu = 0.6$}
    \label{sol12}
  \end{subfigure}
  \begin{subfigure}{5cm}
    \centering\includegraphics[width=5cm]{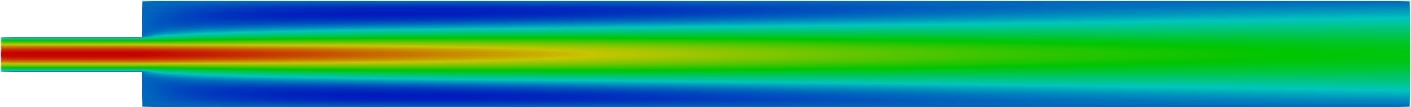}
    \caption{$\nu = 0.3$}
    \label{sol13}
  \end{subfigure}
  
  \begin{subfigure}{5cm}
    \centering\includegraphics[width=5cm]{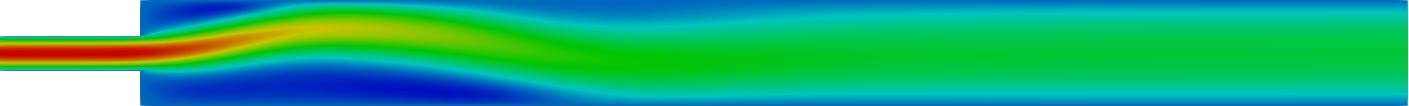}
    \caption{$\nu = 0.6$}
    \label{sol2}
  \end{subfigure}
  \begin{subfigure}{5cm}
    \centering\includegraphics[width=5cm]{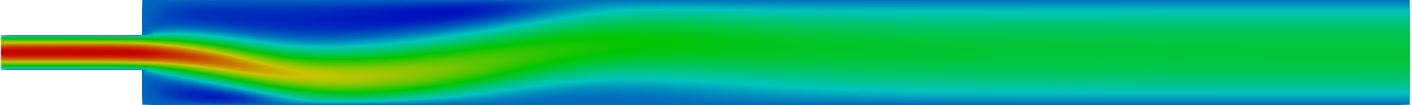}
    \caption{$\nu = 0.6$}
    \label{sol3}
  \end{subfigure}
  \begin{subfigure}{5cm}
    \centering\includegraphics[width=5cm]{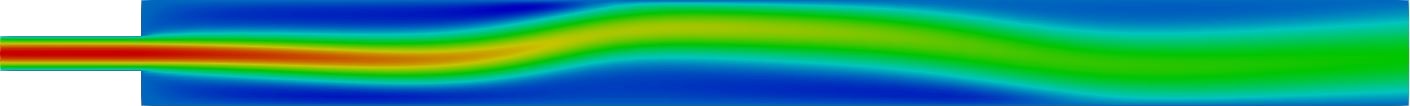}
    \caption{$\nu = 0.3$}
    \label{sol4}
  \end{subfigure}
  
  \begin{subfigure}{5cm}
    \centering\includegraphics[width=5cm]{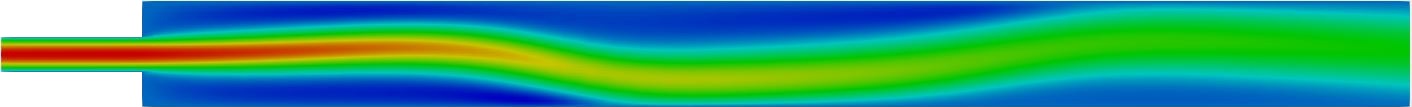}
    \subcaption{$\nu = 0.3$}
    \label{sol5}
  \end{subfigure}
  \begin{subfigure}{5cm}
    \centering\includegraphics[width=5cm]{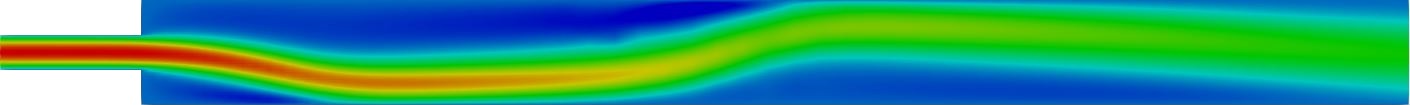}
    \caption{$\nu = 0.3$}
    \label{sol6}
  \end{subfigure}
  \begin{subfigure}{5cm}
    \centering\includegraphics[width=5cm]{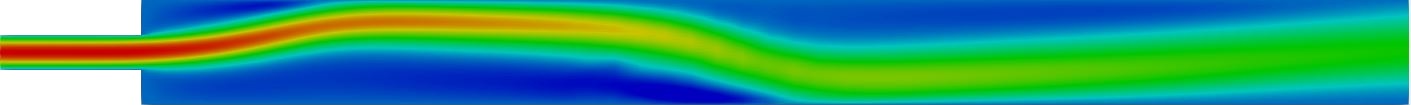}
    \caption{$\nu = 0.3$}
    \label{sol7}
  \end{subfigure}
  
   \begin{subfigure}{5cm}
    \centering\includegraphics[width=5cm]{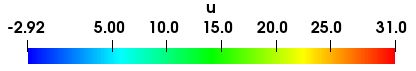}
  \end{subfigure}
  \caption{Nine of the most representative solutions that have been obtained in this work. The colour gradient is associated with the streamwise velocity. Solutions \protect{\subref{sol11}}, \protect{\subref{sol12}} and \protect{\subref{sol13}} belong to branch 1 of figure \ref{1d_diagr}, \protect{\subref{sol2}} and \protect{\subref{sol7}} to branch 2, \protect{\subref{sol3}} and \protect{\subref{sol6}} to branch 3, \protect{\subref{sol4}} to branch 4 and, finally, \protect{\subref{sol5}} to branch 5}
  \label{possible_solutions}
\end{figure}

\subsection{Results with a single varying parameter}
\label{res_1_param}
As previously written, in this section we will discuss the results that can be obtained with only a single parameter. Therefore, we will not consider the efficiency of the computation but we will focus on the application of the described techniques to accurately discretize a bifurcation diagram. In fact, the bifurcation diagrams computed offline and online (that will be described later and are shown in figures \ref{1d_diagr} and \ref{online_1d_diagram}) represent almost the same information. Then, when interested in such a diagram, could simply use the one obtained offline without performing the POD and the online phase. We thus only prove, in this section, that the proposed method is stable and that the deflated continuation method can be used also in a reduced framework without major changes. On the other hand, in section \ref{res_2_param} we will construct bifurcation diagrams with more parameters to show how to exploit the RB method to ensure the efficiency.\\
The first diagrams that we want to show represent the decay of the eigenvalues of the correlation matrix obtained in the POD method. 
In figure \protect{\subref{svdecay1}} one can observe the decay associated with 24 snapshots distributed on three different branches near a bifurcation point, while the entire diagram that contains such branches can be seen in figure \ref{1d_diagr}. Such snapshots have been computed with $\nu\in[0.85,1]$, $s=1$ and with the additional constraint $|\nu_{i+1}-\nu_i|=\Delta\nu_i<\Delta\nu_{max}=0.02\cdot\nu$ to ensure that the approximations obtained with the continuation method are accurate enough. In the following we will refer to this particular choice for the number of snapshots, the viscosity range, the value of $s$ and this particular constraint on $\nu$ as the reference setting. It can be observed that the decay is exponential even if there is a singular point (for instance, similar behaviours have been proved in \cite{exp_decay} and in \cite{kolmogorov} without considering singular points). Moreover, it should be noted that the decay is very similar when the continuation steps are smaller and, therefore, all the snapshots are closer to the bifurcation point (see figure \protect{\subref{svdecay2}}, obtained with $\nu\in[0.91,1]$ and $\Delta\nu_{max}=0.01\cdot\nu$). On the other hand, the decay remains exponential but it is faster or slightly different if, respectively, the snapshots belong to the same branch or if we increase their number from 24 to 100. Note that discarding two branches or increasing the number of snapshots implies that, with the same step sizes, the viscosity varies in a wider range, approximately [0.5,1] in both cases. Such results are important because an exponential decay of the eigenvalues implies, thanks to relation \eqref{pod_error1}, that the approximation error of the reduced spaces exponentially decreases with respect to its own dimension.\\
\begin{figure}
\Large
\centering  
  \begin{subfigure}{5.8cm}
    \centering\includegraphics[width=5.75cm]{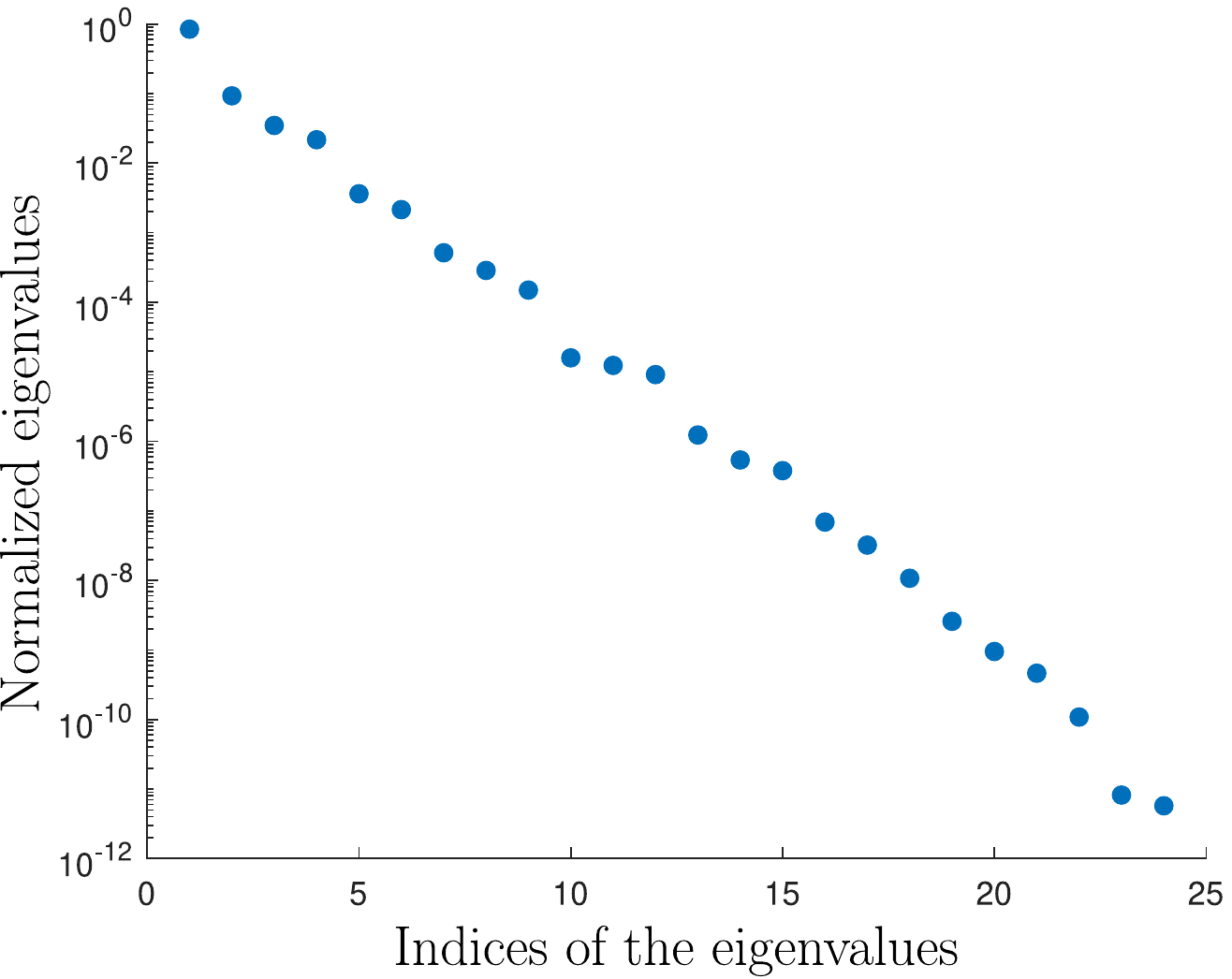}
    \caption{Reference setting}
    \label{svdecay1}
  \end{subfigure}
  \begin{subfigure}{5.8cm}
    \centering\includegraphics[width=5.75cm]{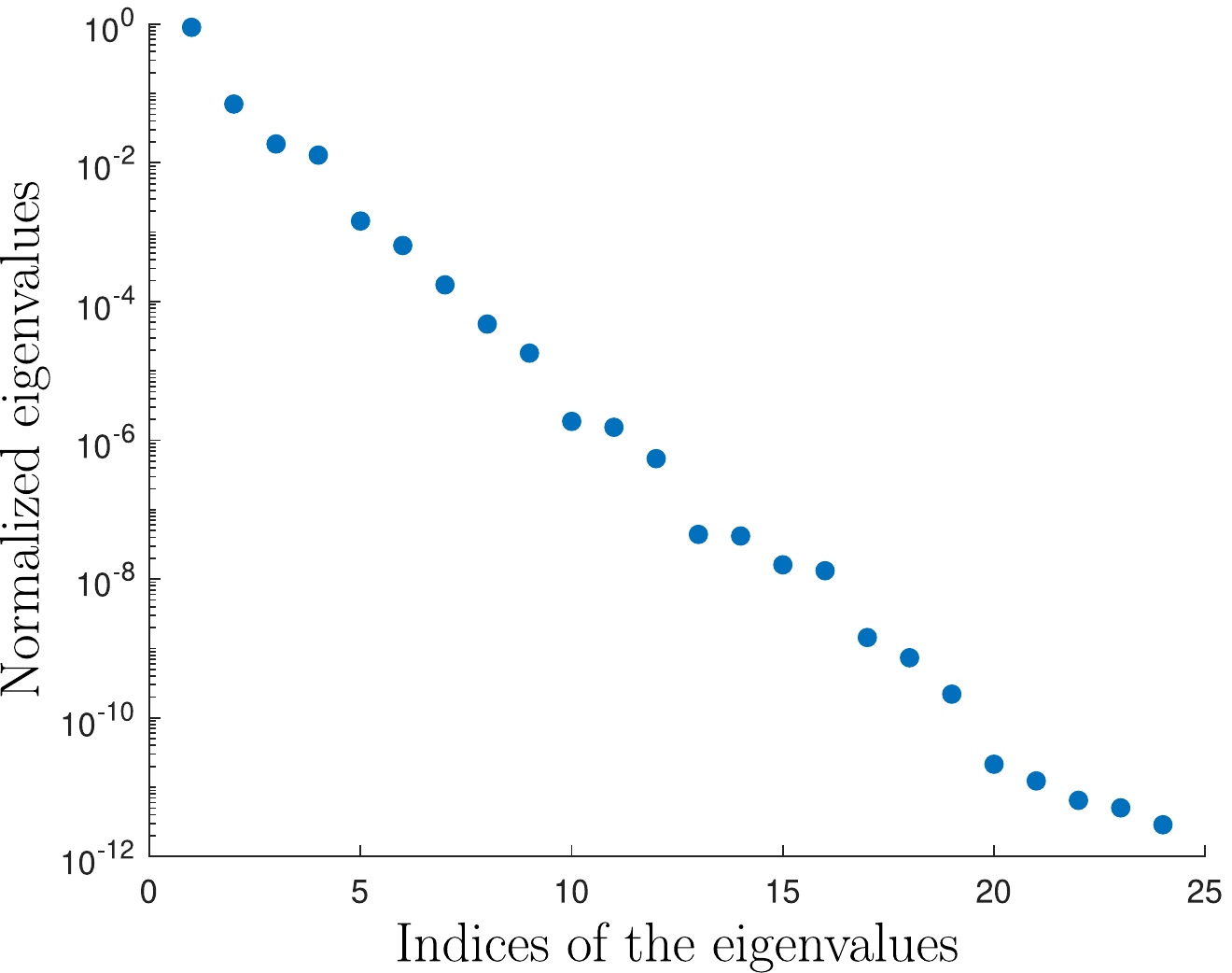}
    \caption{Reference setting computed with $\Delta\nu_{max}=0.01\cdot\nu$}
    \label{svdecay2}
  \end{subfigure}
  
  \begin{subfigure}{5.8cm}
    \centering\includegraphics[width=5.75cm]{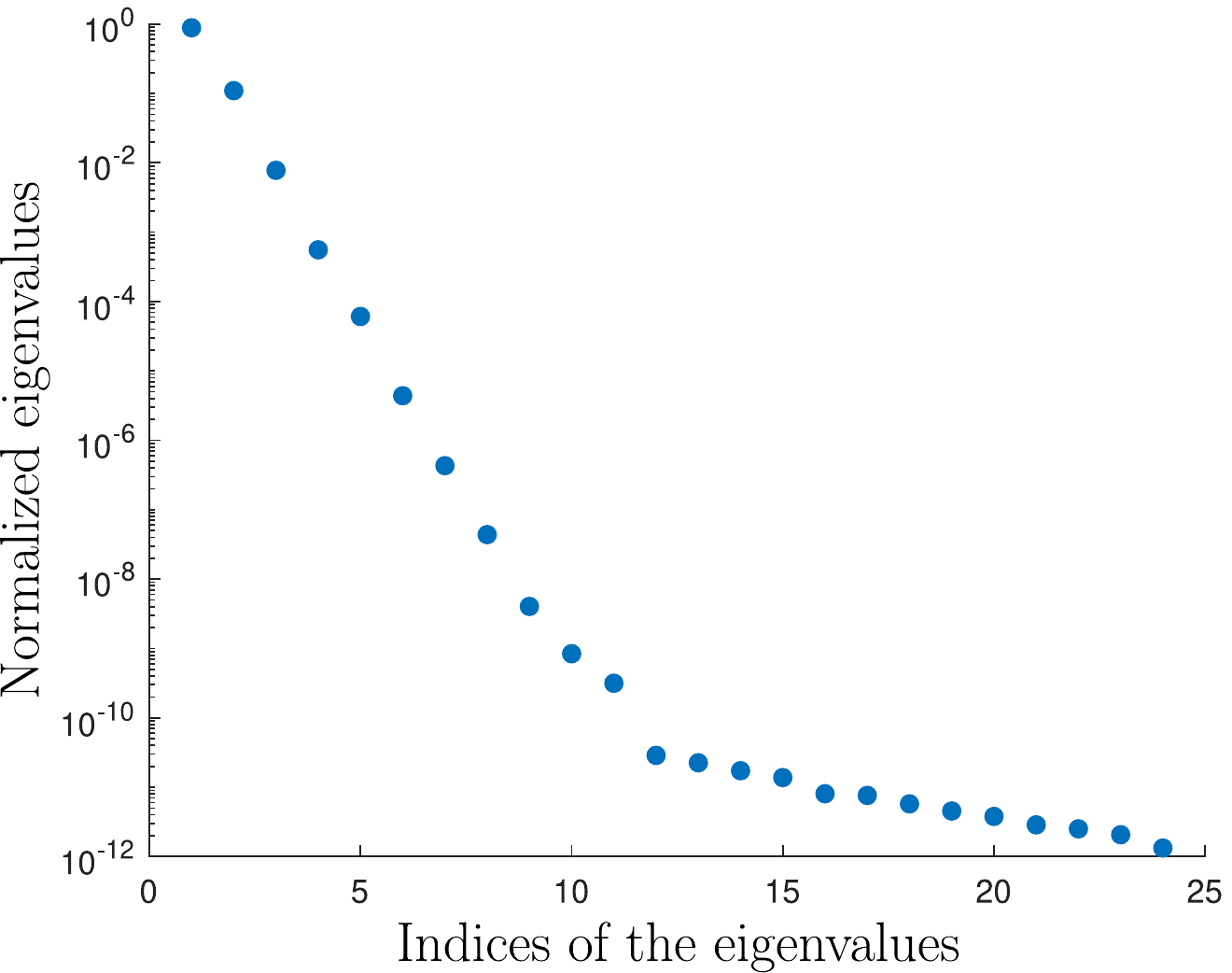}
    \caption{Reference setting with all the snapshots belonging to the same branch}
    \label{svdecay3}
  \end{subfigure}
  \begin{subfigure}{5.8cm}
    \centering\includegraphics[width=5.75cm]{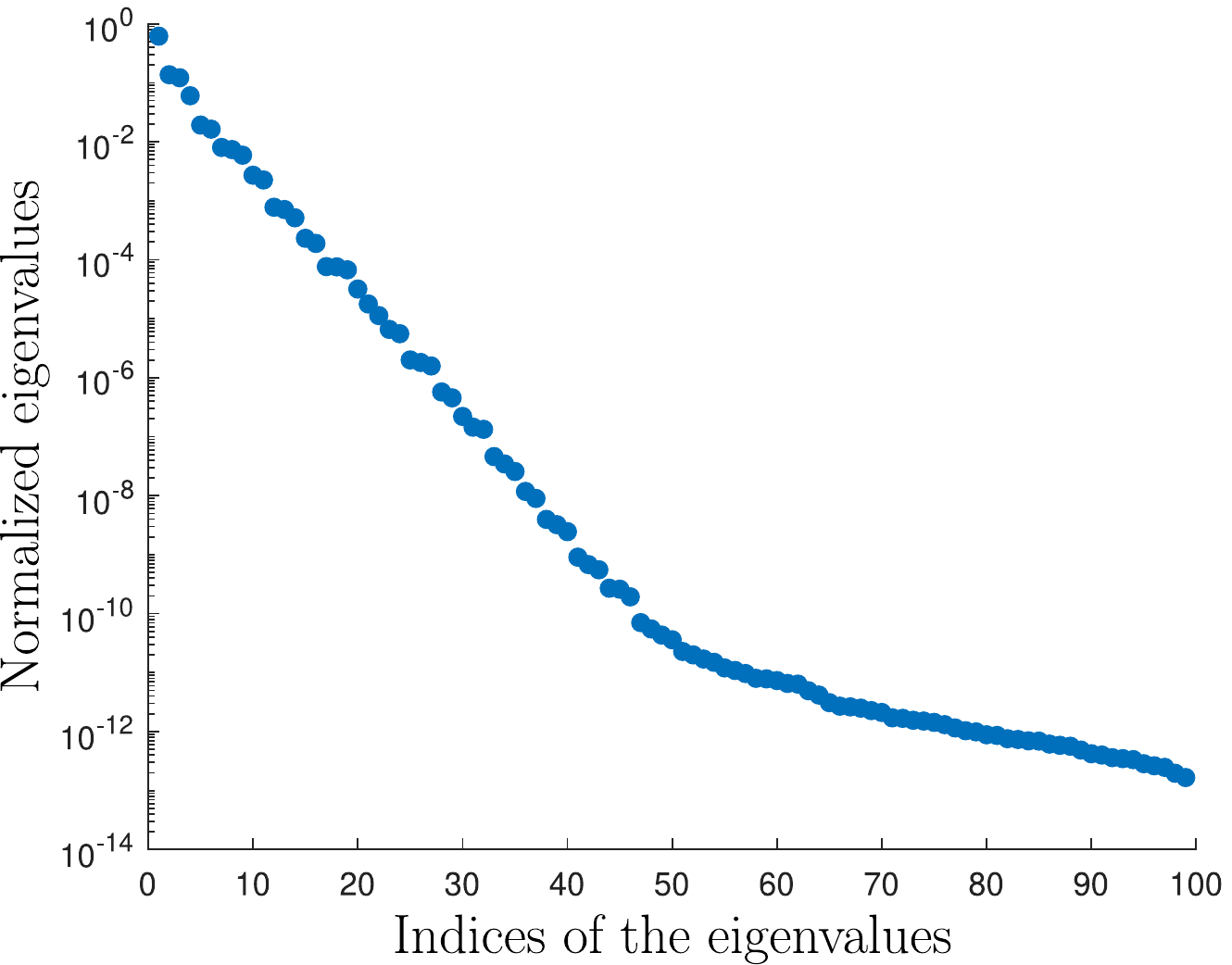}
    \caption{Reference setting with 100 snapshots}
    \label{svdecay4}
  \end{subfigure}
  \caption{Decays of the eigenvalues of the correlation matrix used in the POD method. The reference setting consists in 24 snapshots belonging to the first three branches of the diagram in figure \ref{1d_diagr} (higher values of $\nu$) computed with $\nu\in[0.85,1]$, $s=1$ and $\Delta\nu_i<\Delta\nu_{max}=0.02\cdot\nu$}
  \label{eig_decay_1d}
\end{figure}The next figure shows the entire one-dimensional bifurcation diagram which we will discuss (figure \ref{1d_diagr}). It is computed during the offline phase and each point corresponds to a snapshot.
\begin{figure}[h]
\centering
\includegraphics[height=10cm, width=10cm, angle=0, keepaspectratio]{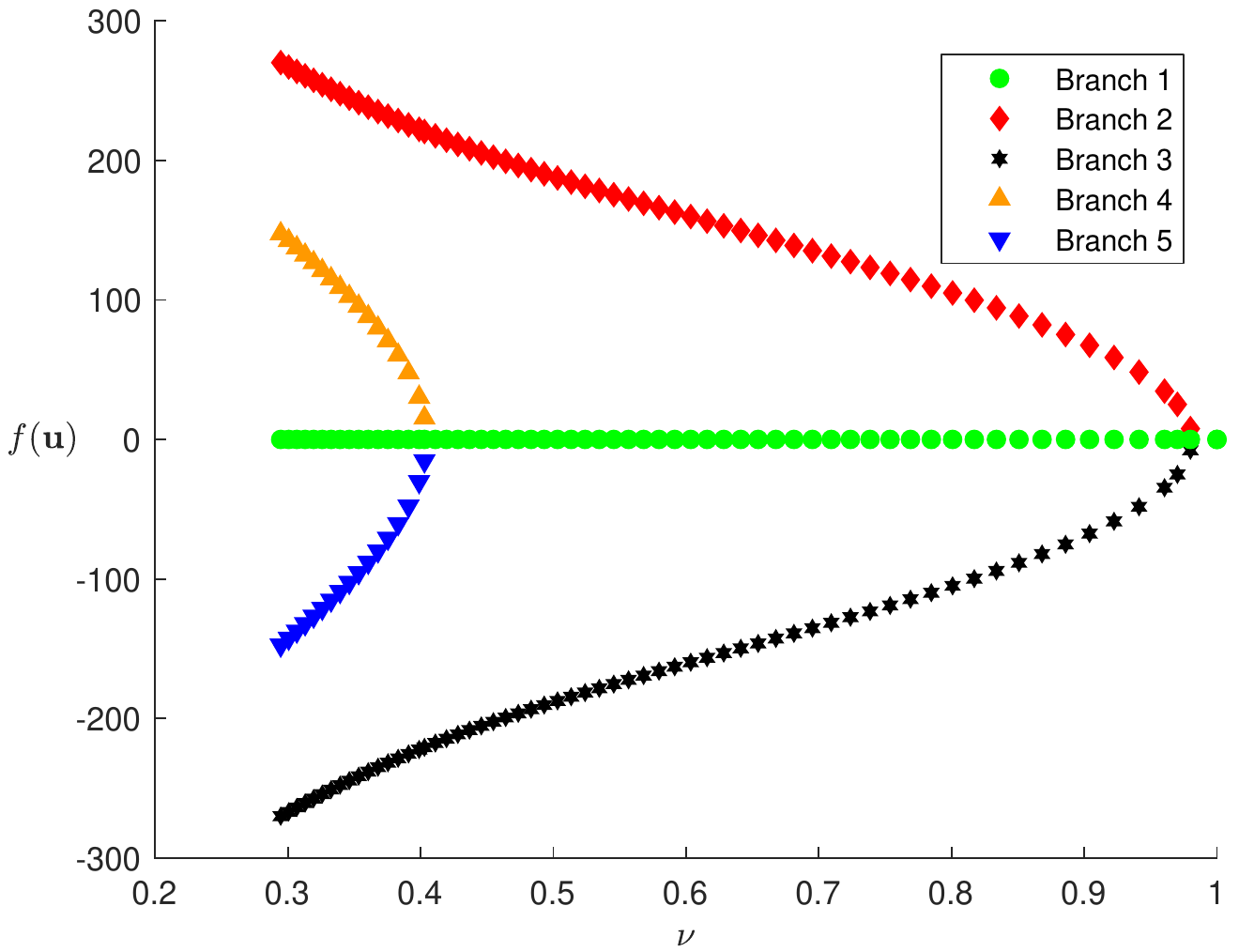}
\caption{Bifurcation diagram computed during the offline phase, each point is associated with a snapshot and each branch is characterized by a different colour and a different marker}
\label{1d_diagr}
\end{figure}It can be observed that it includes two bifurcation points and five different branches. In order to compute it, we decided to use as first solution the one obtained for $\nu=1$ and, then, we decreased the viscosity computing the solutions with the deflated continuation method. Note that the pseudo-arclength continuation automatically selects the best $\Delta\nu_i$ values for each branch, but we imposed the additional constraint that the solutions on different branches have to be associated with the same viscosities to make the deflation more effective and that, in order to avoid too wide steps, $\Delta\nu_i<\Delta\nu_{max}=0.02\cdot\nu$. We also highlight that, if one is not interested in the entire diagram but only in the position of the bifurcation points, it is convenient to perform an eigenvalues analysis as discussed in \cite{pichi_bif} and in \cite{eigen_study}.\\
Since, as discussed in the introduction of this section, the output functional $f(\cdot)$ can be considered as a measure of the asymmetry of the solutions, one can observe that the solutions in figures \protect{\subref{sol11}}, \protect{\subref{sol12}} and \protect{\subref{sol13}} are associated with the middle branch. Moreover, the ones in figures \protect{\subref{sol2}} and \protect{\subref{sol7}} are two characteristic velocity fields of the first and of the second part of the upper branch that is born from the first bifurcation point (the one on the right, since we decrease the viscosity). \\
Consequently, the solution in figures \protect{\subref{sol3}} and \protect{\subref{sol6}} are associated with the lower branch because they are the mirrored solutions of the latter and, finally, the fields in figures \protect{\subref{sol4}} and \protect{\subref{sol5}} are representative of the solutions on the last two branches. Obtaining such a diagram is very expensive because many different solutions have to be computed with the full order solver. In fact it has been obtained using 224 snapshots with $N^\delta =$ 7372 degrees of freedom and, to obtain them, each step of the iterative solver spent almost 0.67 seconds. We remark that with 7372 degrees of freedom we were able to compute accurate solutions thanks to the SEM. However, if one uses other discretization techniques or requires a very high accuracy, the computational cost can significantly increase. For instance, in \cite{local_rom}, a very similar model has been computed with the FEM (in particular with the Taylor-Hood elements) exploiting 90876 degrees of freedom to obtain the desired accuracy. Note that the 7372 degrees of freedom are associated with the mesh shown in figure \ref{mesh} and with ansatz functions of order 12. Therefore the $L^2$ error is of order $\mathcal O(h^p)$ where $p=12$ and $h$ can be computed as the ratio between the diameter of the biggest element and the domain diameter, such a ratio is approximately 0.05. \\
It is important to note that, even if the computation of this diagram is very expensive, such a result can be obtained during the online phase much more efficiently. The obtained diagram is shown in figure \ref{online_1d_diagram}, it can be seen that it is possible to reconstruct both the bifurcation points and all the branches. Such a figure has been obtained using the simple continuation and dividing by 20 the step sizes used during the offline phase. Moreover, since in this simulation the solutions were associated with $N^{rb} =$ 37 degrees of freedom, the iterative solver approximately spent only $10^{-5}$ seconds for every iteration with a speedup of about 4 order of magnitude. The statistics in table \ref{tab_time1} can be used to observe that, if a bifurcation diagram is discretized with at least 770 solutions, then the described approach is more efficient than the direct computation of a diagram without the RB method. Finally, we remark that these quantities have been obtained exploiting a prior knowledge obtained from previous experiments about the position of the bifurcation points to use the deflation only when new branches can be found. However, without such a knowledge one should use the deflation at each step of the continuation, extremely increasing the computational cost of the process and making the described approach even more convenient.\\
Unfortunately, because of the data structures used in \textit{Nektar++} to perform static condensation, at the moment some of the computations of the online phase depended on $N^\delta$ and, therefore, the time required to compute, on average, a solution given the previous ones significantly increased ruining the speedup obtained in a single iteration. In fact the online solver could generally converge in 4 iterations, but the average time obtained dividing the total required time needed to get the complete diagram by the number of computed solutions is much higher than 4 times the time required to perform a single step of the iterative solver. Approaches to obtain a full-decoupling from $N^\delta$, exploiting a reduced change of basis matrix, will be analyzed in future works.``
\begin{table}[h] \centering \fontsize{8.5}{9.9} \selectfont
\begin{tabular}{c c c c c c c} 
\hline \\[-0.2cm]
&$T_d$&$N_d$ &$T_{1s}$&$T_{POD}$&$T_{dN}(N)$&$T_{1i}$\\ [0.5ex]
\hline \\[-0.2cm]
Offline & 1292 s & 224 & 5.77 s & 2995 s & $(5.77\cdot N)$ s & 0.67 s\\ 
Online & 293.85 s & 1492 & 0.20 s & / & $(4287+0.20\cdot N)$ s &  $3.21\cdot10^{-5}$ s\\[1ex] 
\hline\\[-0.8ex] 
\end{tabular}
\caption{Computational costs of the offline and online phases when a single parameter is involved. These quantities have been obtained exploiting a previous knowledge on the position of the bifurcation points. $T_d$ is the time required to compute the entire diagram (with $N_d$ solutions), $T_{1s}$ and $T_{1i}$ are the average times required to compute a single solution and to perform a single iteration of the iterative solver, while $T_{POD}$ is the time required by the POD and $T_{dN}(N)$ the one needed to compute an arbitrary diagram with $N$ solutions. In the online phase $T_{dN}(N)$ is computed summing the time required to compute the snapshots, to perform the POD and to compute $N$ reduced solutions. The described approach, with this offline phase, is more efficient than the classical one when the bifurcation diagram is discretized with $N\ge770$ solutions because $N=770$ is the first value such that $T_{dN}(N)$ of the offline phase is higher than its online counterpart.}
\label{tab_time1}
\end{table}\\In order to quantify the accuracy of the obtained result, we decided to perform an empirical error analysis. Thus we reprojected the reduced solutions on the full order space and used them as initial guesses for the iterative solver. This way we compared the obtained full order solutions with the reduced ones and we computed the associated relative error that is shown in figure \ref{1d_error}. In order to properly interpret the diagram, it is important to remark that the tolerance used by the offline iterative solver was $10^{-6}$ and, therefore, the relative error is almost always very close to such a quantity. However, coherently with the fact that one of the hypothesis of the reduced basis method regards the presence of a smooth solution manifold, the error increases in the neighbourhoods of the bifurcation points. This phenomenon is in accordance with the fact that it is more complex to converge to solutions very close to singular points without good initial guesses (that are never available when one wants to discover unknown branches). It can also be observed in figure \ref{online_1d_diagram}, where the asymmetrical branches were not always immediately found.\\
\begin{figure}[h]
\centering
\includegraphics[height=10cm, width=10cm, angle=0, keepaspectratio]{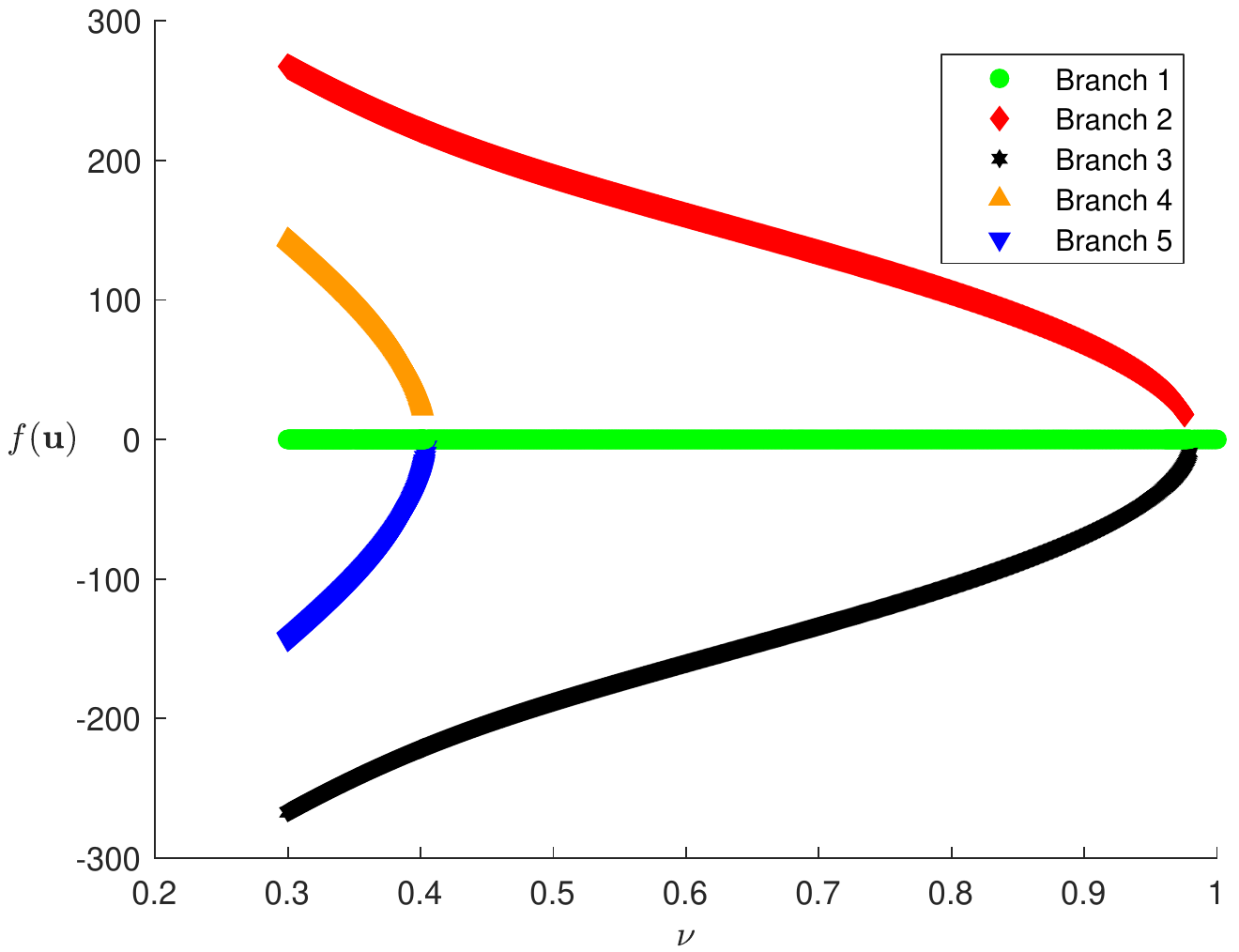}
\caption{Bifurcation diagram efficiently reconstructed during the online phase}
\label{online_1d_diagram}
\end{figure}Moreover, it is possible to observe the behaviours of the average and maximum relative errors in table \ref{tab1d}. The first column is associated with a neighbourhood of the bifurcation point near $\nu\approx0.4$, the second one to the portion of the domain between the singular points, the third to a neighbourhood of the other bifurcation point and, finally, the fourth to the entire diagram. It can be noted that, even if the bifurcation points strongly affect the accuracy, the global average error is only one order of magnitude above the solver threshold.
\begin{figure}
\centering  
  \begin{subfigure}{7.0cm}
    \centering
    \includegraphics[width=6.85cm]{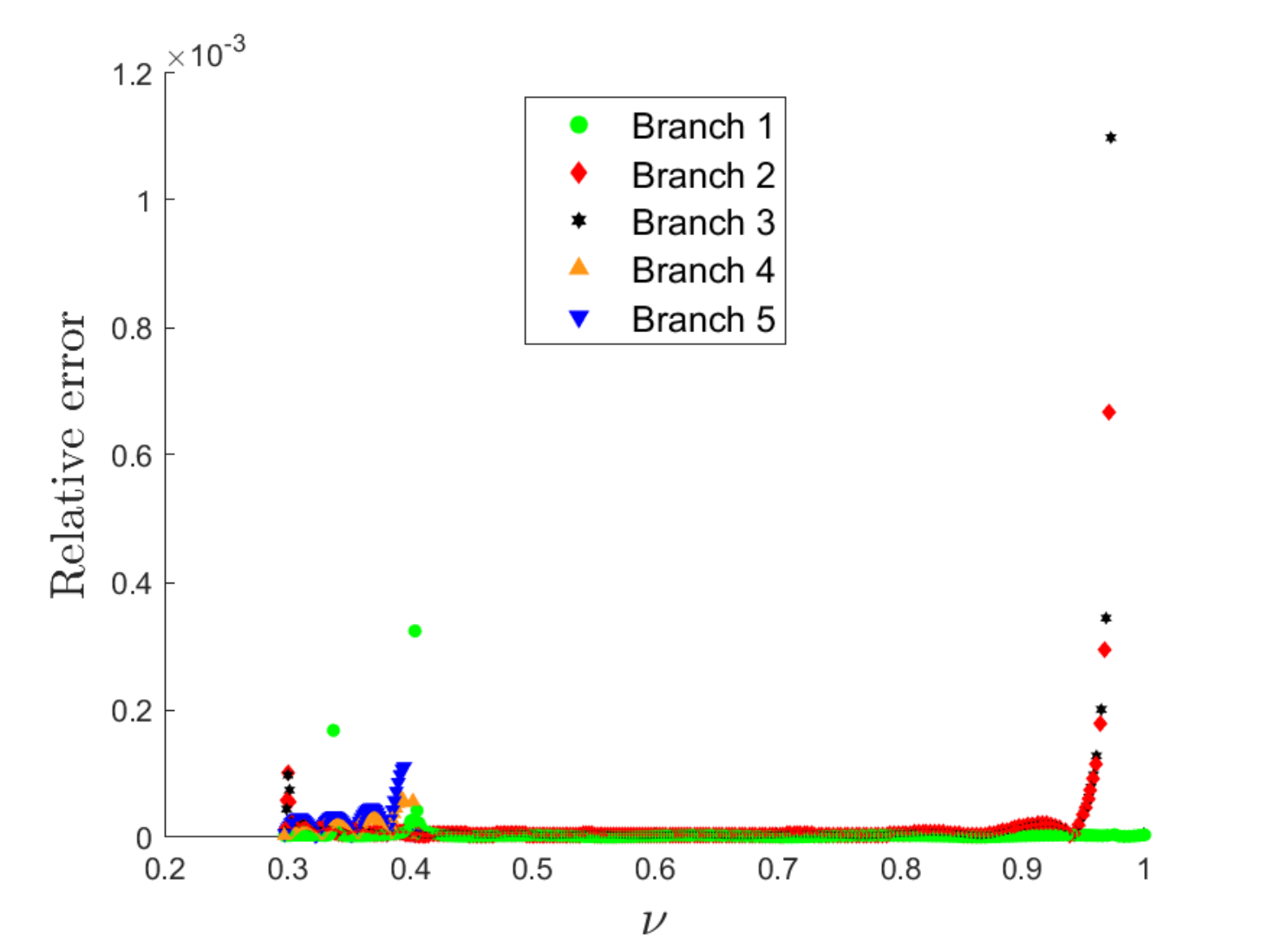}
    \caption{Uniform scale}
    \label{error_1d_unif}
  \end{subfigure}
  \begin{subfigure}{7.0cm}
    \centering
    \includegraphics[width=6.85cm]{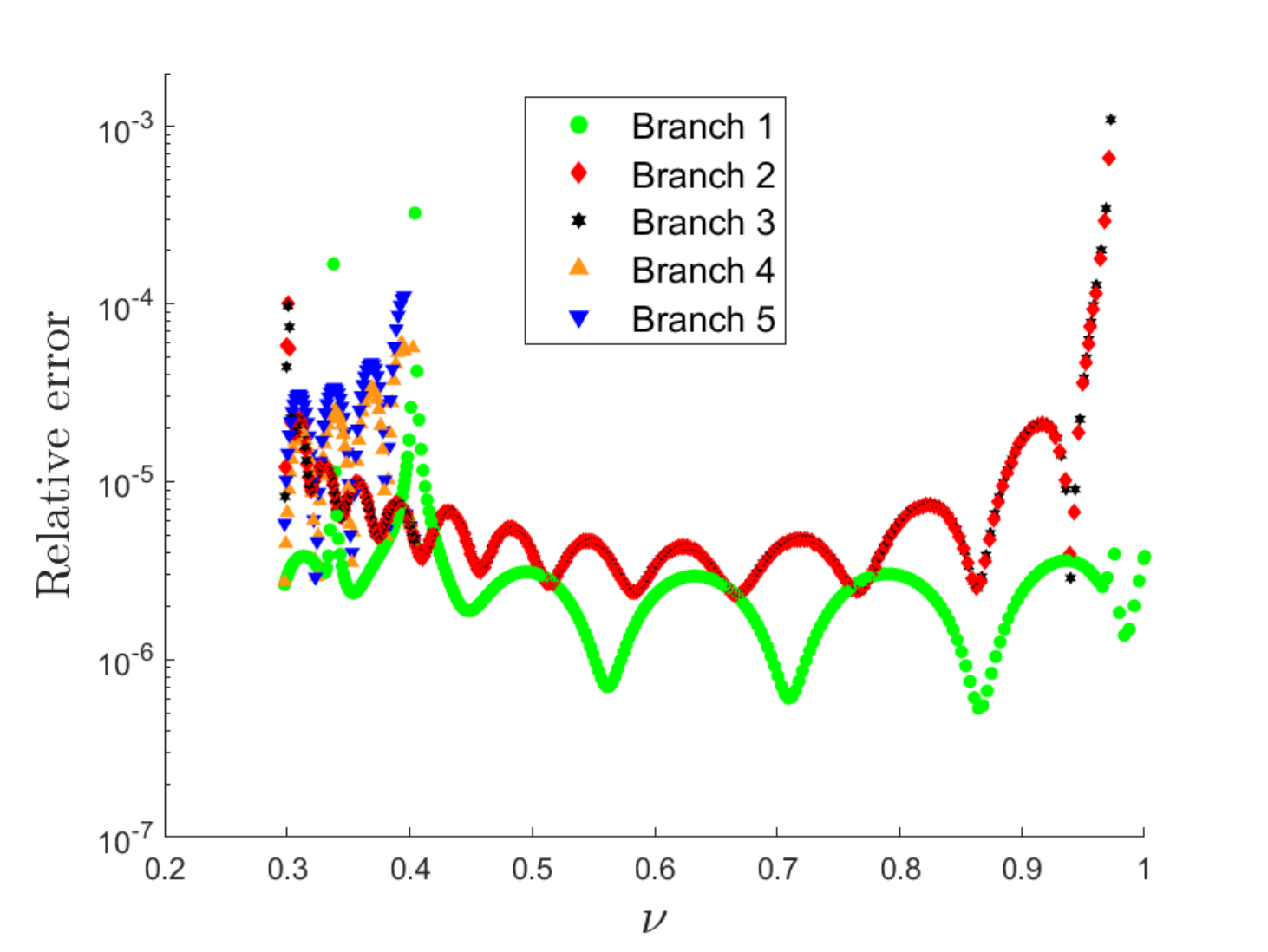}
    \caption{Logarithmic scale}
    \label{error_1d_log}
  \end{subfigure}
  \caption{Relative error between reduced and full order solutions of figures \ref{1d_diagr} and \ref{online_1d_diagram}} 
  \label{1d_error}
\end{figure} \begin{table}[h]
\centering 
\begin{tabular}{c c c c c} 
\hline \\[-0.2cm]
Quantity&\hspace{0.15cm}$\nu\in\Omega_0^\nu$\hspace{0.15cm}&\hspace{0.15cm}$\nu\in\Omega_1^\nu$\hspace{0.15cm}&\hspace{0.15cm}$\nu\in\Omega_2^\nu$\hspace{0.15cm}&\hspace{0.15cm}$\nu\in P$\hspace{0.15cm}\\ [0.5ex]
\hline \\[-0.2cm]
Average error & $1.53\cdot10^{-5}$ & $3.56\cdot10^{-6}$ &$ 5.45\cdot10^{-5}$ &  $1.24\cdot10^{-5}$\\ 
Maximum error &$3.24\cdot10^{-4}$& $1.67\cdot10^{-5}$ & $1.10\cdot10^{-3}$ &  $1.10\cdot10^{-3}$\\[1ex] 
\hline \\[-0.8ex]
\end{tabular}
\caption{Average and maximum relative error of the solutions of diagram \ref{online_1d_diagram} with respect to the corresponding full order solutions. The following notation has been used: $\Omega_0^\nu=[0.3,0.45]$, $\Omega_1^\nu=(0.45,0.9]$, $\Omega_2^\nu=(0.9,1]$, $P=\Omega_0^\nu\cup\Omega_1^\nu\cup\Omega_2^\nu=[0.3,1]$.} 
\label{tab1d}
\end{table}\\It is interesting to note that the error exponentially decreases with respect to the dimension of the reduced space, both in smooth regions and close to the singular points. This phenomenon can be observed in figure \ref{error_over_rbsize_mm}, where the error is associated with solutions in a small neighbourhood of the first bifurcation point, i.e. a region where the hypothesis of smooth solution manifold required by the RB method does not hold. To obtain such a decay we computed 30 snapshots with $\nu\in[0.825,1]$, $s=1$ and $\Delta\nu_{max}=0.02\cdot\nu$, we then generated 30 different reduced spaces and compared the online solutions with the corresponding full order ones. It is also important to observe that the error, computed over 300 reduced solutions, stops decreasing when it is close to $10^{-10}$ because such a quantity has been used as the tolerance in the offline and online iterative solvers. The fact that the maximum error decreases almost as the average one implies that the accuracy of the less accurate solution (the closest one to the singular point) exhibits the same behaviour as the other ones.

\begin{figure}
\centering
\includegraphics[height=9.6cm, width=9.6cm, angle=0, keepaspectratio]{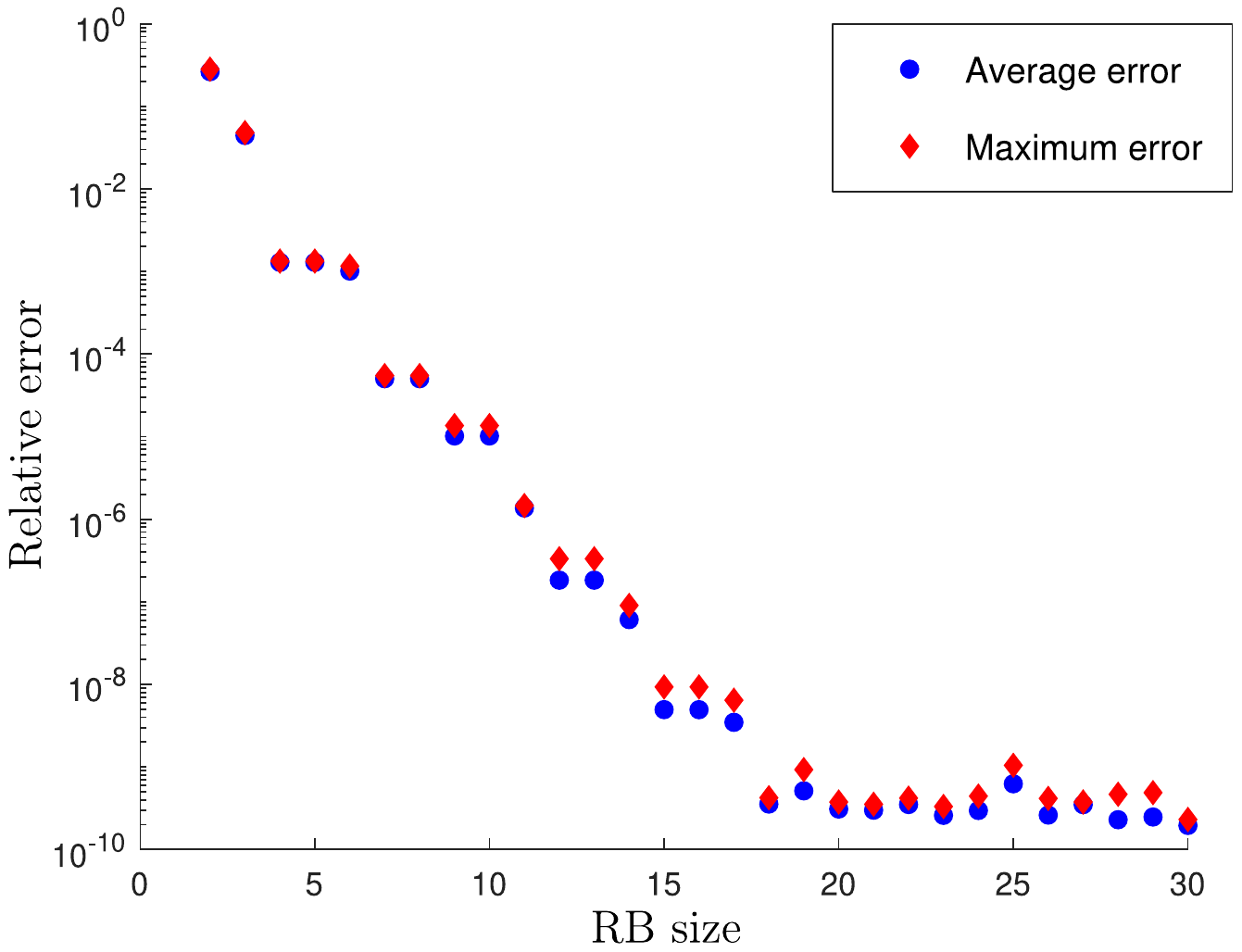}
\caption{Relative errors over the dimension of the reduced space }
\label{error_over_rbsize_mm}
\end{figure}

\subsection{Results with multiple parameters}
\label{res_2_param}
Finally, in this section we will address the problem of efficiency. In the previous section we showed that the eigenvalues decay of the POD matrix is the expected one and that, exploiting the continuation and the deflation method, it is possible to construct a bifurcation diagram in the offline phase and to reconstruct it accurately and efficiently during the online one. \\
In section \ref{res_1_param}, the diagram discretized during the online phase loses part of its importance because the same branches have to be discretized also in the offline one. Otherwise, if one wants to obtain it only in the online one without exploiting the continuation and the deflation offline, the reduced space would be too small and the secondary branches could not be found online. For instance, let us assume that during the offline phase one computes all the symmetric solutions, i.e. the ones such that $f(\uu) = 0$, then the reduced space is able to generate only symmetric solutions and, therefore, the other branches can not be discovered.\\
Furthermore, if one is interested in computing bifurcation diagrams with more than one parameter, the number of solutions required to discretize it would significantly increase and the computational cost of the high-order simulation would become prohibitive. We remark that the coupling of the aforementioned techniques allows us to efficiently reconstruct a representation of the $N$-dimensional manifold induced by the solutions, that would be infeasible without reduction strategies. Therefore, in order to compute a diagram letting vary more parameters, we decided to slightly change the described approach, computing only few one-dimensional bifurcation diagrams during the offline phase (in this work we computed only $N_{\text{off\_diag}}=2$ or $N_{\text{off\_diag}}=3$ offline diagrams), and refining the grid associated with the second parameter only in the online phase. Such an approach can be generalized, when even more parameters are involved, by computing a small set of one-dimensional diagrams during the offline phase and generating all the other dimensions of the complete diagram only online.\\
As previously written, the second parameter in this work is a scaling $s$ of the inlet Dirichlet boundary condition. The bifurcation diagram that can be obtained with the described approach is shown in figure \ref{2d_diagr}. It can be observed that it is possible to entirely reconstruct the two-dimensionality of the diagram with both the bifurcation points, whose position changes according to the value of $s$. The obtained diagram is coherent with the one in figures \ref{1d_diagr} and \ref{online_1d_diagram} and with the nature of the two involved parameters. In fact the inlet boundary conditions is strongly related to the Reynolds number ($Re=UL/\nu$) and, since such a non-dimensional quantity is the one responsible for the bifurcation, one could expect that the bifurcation points are moved along straight lines in the $s$-direction because multiplying $\nu$ by a scalar $c$ is equivalent, in terms of $Re$, to dividing $s$ by $c$. Moreover, this property explains the fact that, when $s$ increases from 0.8 to 1, the bifurcation phenomena are anticipated, i.e. the bifurcation points are associated with higher viscosities.\\
\begin{figure}
\centering
\includegraphics[height=11.8cm, width=11.8cm, angle=0, keepaspectratio]{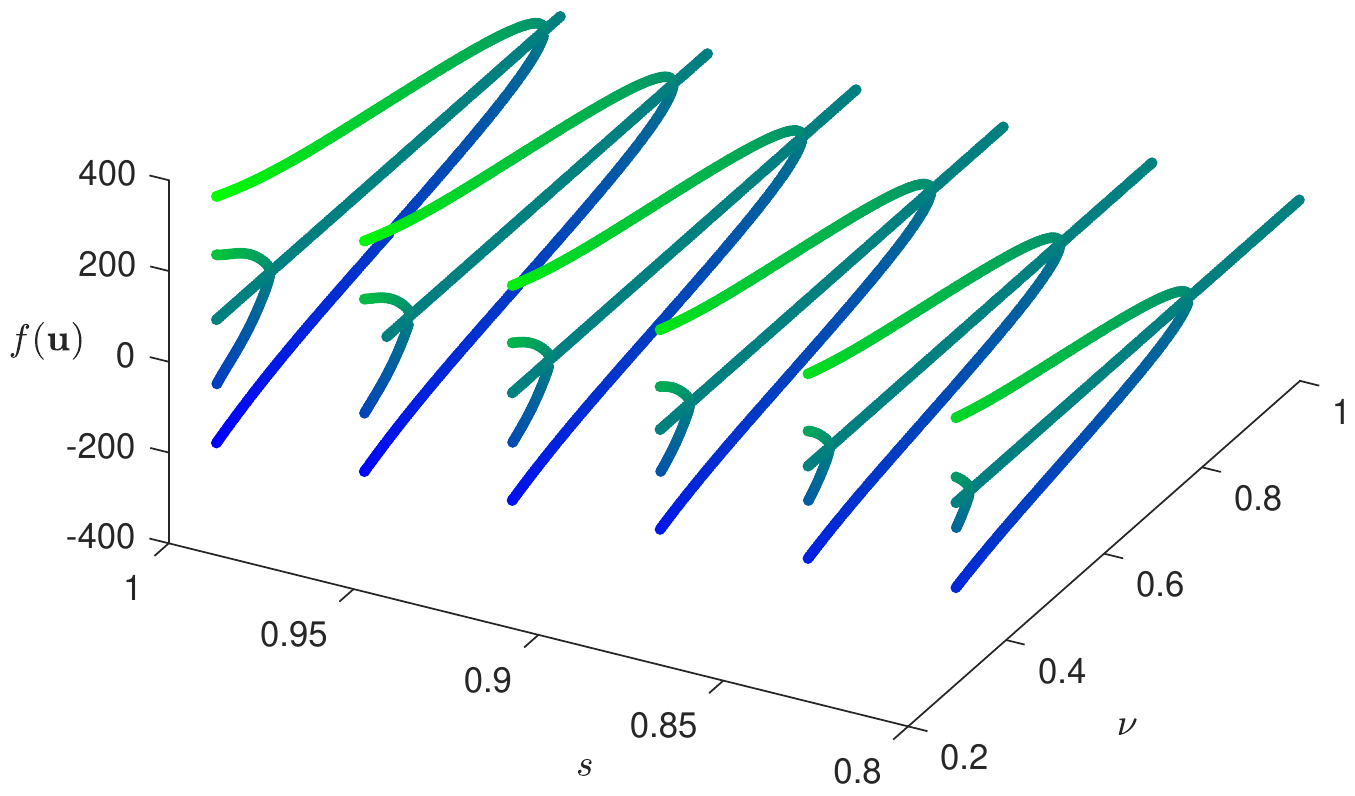}
\caption{Bifurcation diagram efficiently obtained in the online phase with two parameters and two bifurcations. The colour gradient remarks the value $f(\uu)$ in each point in order to allow the reader to more easily interpret the diagram}
\label{2d_diagr}
\end{figure}We wanted to consider a second parameter that could add new information to the one-dimensional diagram but, since these two parameters are so deeply related to $Re$, we decided to analyze the POD eigenvalues decay to observe if a behaviour different than the one observed in the one-dimensional case could be found. A significantly stronger decay would have meant that $\nu$ and $s$ contained the same information and, therefore, the same solutions that can be obtained varying the viscosity could be obtained properly rescaling the ones obtained varying $s$. \\
We analyzed three different scenarios. Firstly, we considered snapshots belonging to a diagram with 2 values of $s$ and 15 of $\nu$, that is very similar to the actual offline phase and it is reported in figure \protect{\subref{svdecay5}}. In such a setting, the observed decay follows the behaviour of the reference one-dimensional setting. Secondly, in figure \protect{\subref{svdecay6}} we used only 1 value of $\nu$ and 30 of $s$ obtaining 30 solutions on the same branch. Consequently, this decay agrees to the one in figure \protect{\subref{svdecay3}}. Finally, in figure \protect{\subref{svdecay7}} we decided to use a mixed approach considering 5 values of $s$ and 6 of $\nu$, the obtained decay is slightly faster than the first one, because of the relation between the parameters, but it is still exponential.\\
\begin{figure}
\centering  
  \begin{subfigure}{6.6cm}
    \centering\includegraphics[width=6.4cm]{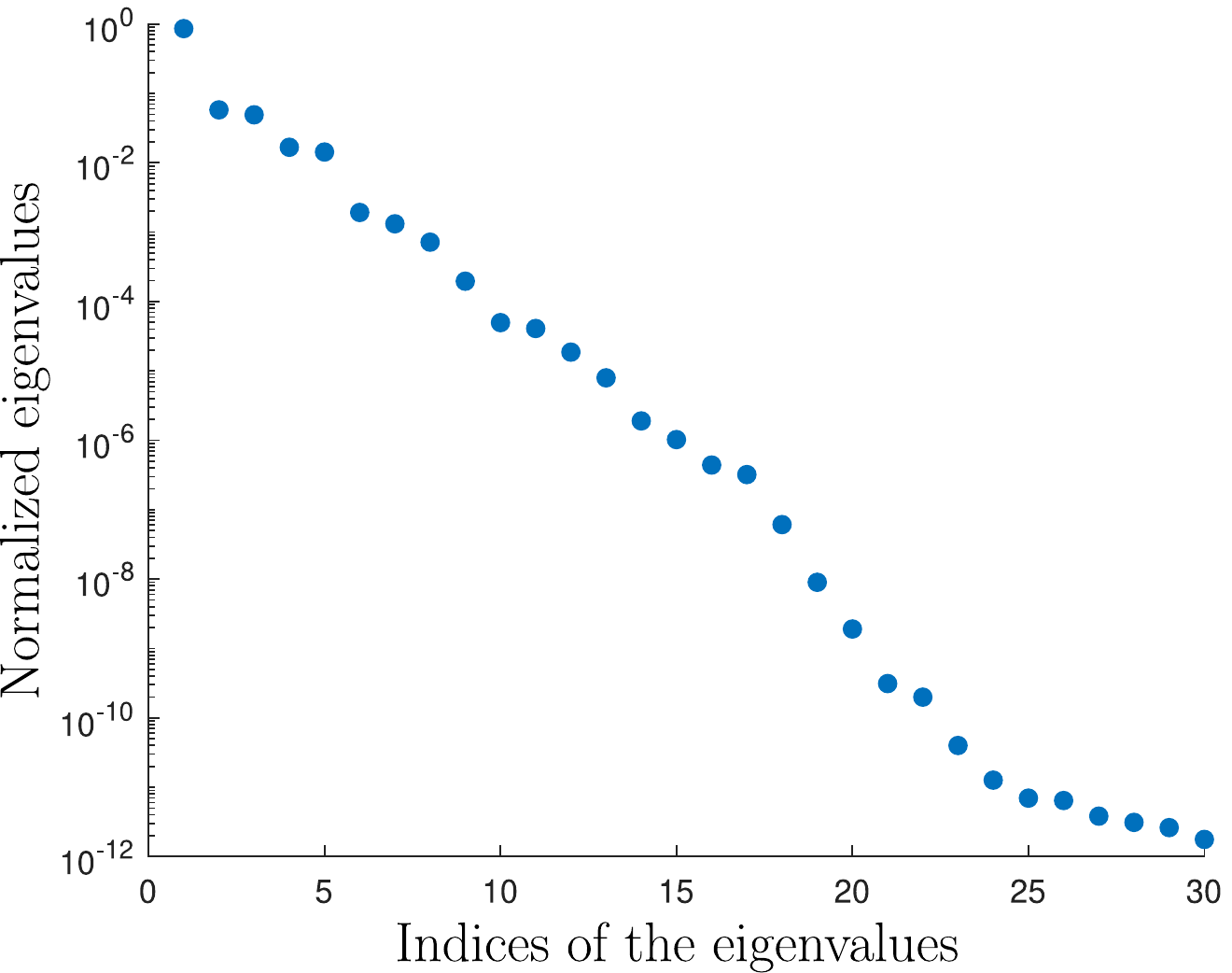}
    \caption{2 values of $s$, 15 values of $\nu$}
    \label{svdecay5}
  \end{subfigure}
  \begin{subfigure}{6.6cm}
    \centering\includegraphics[width=6.4cm]{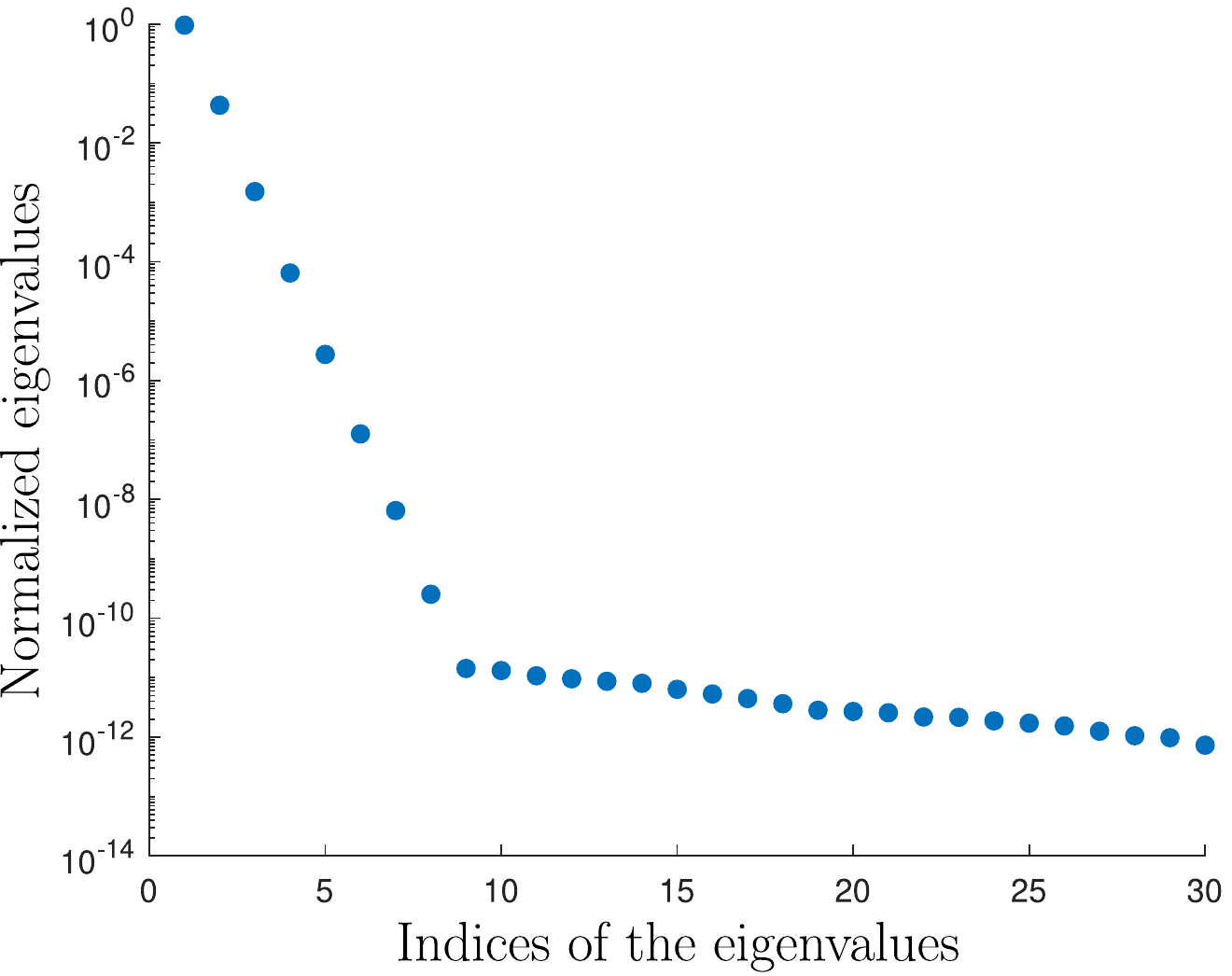}
    \caption{30 values of $s$, 1 value of $\nu$}
    \label{svdecay6}
  \end{subfigure}
  \begin{subfigure}{6.6cm}
    \centering\includegraphics[width=6.4cm]{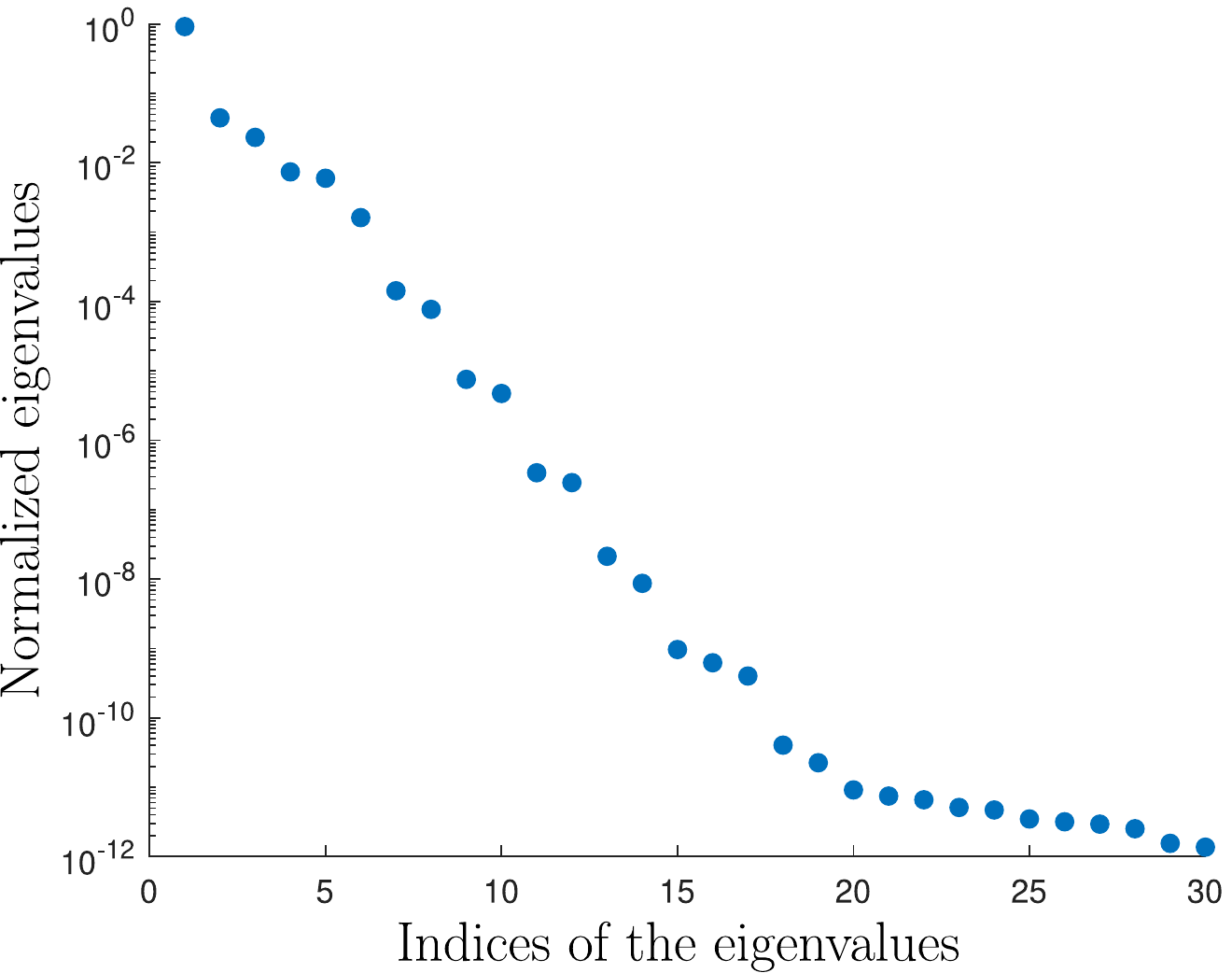}
    \caption{5 values of $s$, 6 values of $\nu$}
    \label{svdecay7}
  \end{subfigure}
  \caption{Eigenvalues decay obtained with different settings of the model obtained using the scaling as a parameter with different sampling for the snapshots}
  \label{eig_decay_2d}
\end{figure}As remarked in the previous section, the online solver is much more efficient. In fact, even if a single step during the offline phase lasted for 0.67 seconds on average, each single step could be performed in approximately $10^{-5}$ seconds during the online one. It is important to observe that there are about four order of magnitude of difference but this quantity does not depend on the number of involved parameters. Subsequently, the described approach is more and more convenient when the number of parameters increases because much more solutions are required to properly discretize the entire bifurcation diagram. For instance, let us consider the diagram in figure \ref{1d_diagr} that is generated by 224 snapshots, while its central branch only contains 64 different solutions. If one wants to compute a two-dimensional bifurcation diagram with approximately the same discretization level on both the dimensions, the required solutions will be about $224\times64$ and, therefore, 14336 solutions are needed. Moreover, if one considers $n$ parameters, the number of required solutions increases to $224\times64^{n-1}$. The cost associated with such a computation is prohibitive even for $n=2$ if a reduced order model is not taken into account. However, it can be drastically reduced with the described approach, in fact, if one considers only two points in each dimension apart from the first one, the number of snapshots decreases to $224\times N_{\text{off\_diag}}^{n-1}$, while the remaining solutions are computed during the online phase. We claim that the computational cost is significantly reduced because the number of required full order solutions is much smaller than the one required to generate a complete $n$-dimensional diagram, because $N_{\text{off\_diag}}$ is very small. On the other hand, in figure \ref{2d_diagr} we decided to refine the grid associated with the viscosity, while fixing only 6 nodes for the scaling in order to properly visualize the diagram, obtaining 16970 different solutions.\\
Once more, in table \ref{tab_time2}, one can observe that, after a very expensive offline phase, the computation of each solution is very efficient during the online one. However, in order to improve the stability of the online solver, we decided to use shorter continuation steps during the offline phase obtaining more solutions (with the same step sizes we would have approximately obtained $2\cdot224$ snapshots because we considered $N_{\text{off\_diag}}=2$), significantly increasing the computational cost of such a phase. Nevertheless, we can observe that the method is efficient if a diagram is discretized with at least 2365 solutions but, in two-dimensional diagrams, this quantity is very low and, therefore, the entire process is very efficient in most cases. Similar considerations as the ones associated with table \ref{tab_time1} hold also for the statistics in table \ref{tab_time2}.
\begin{table}[h] 
\centering 
\fontsize{8.5}{9.9} \selectfont
\begin{tabular}{c c c c c c c} 
\hline \\[-0.2cm]
&$T_d$&$N_d$ &$T_{1s}$&$T_{POD}$&$T_{dN}(N)$&$T_{1i}$\\ [0.5ex]
\hline \\[-0.2cm]
Offline & 5429 s & 937 & 5.79 s & 7764 s & $(5.79\cdot N)$ s & 0.67 s\\ 
Online & 3564 s & 16970 & 0.21 s & / & $(13193+0.21\cdot N)$ s & $5.43\cdot10^{-5}$ s\\[1ex] 
\hline\\[-0.8ex] 
\end{tabular}
\caption{Computational costs of the offline and online phases when two parameters are involved. The notation of table \ref{tab_time1} has been used. The described approach is more efficient than the classical one when the bifurcation diagram is discretized with $N\ge2365$ solutions. These quantities have been obtained exploiting a previous knowledge on the position of the bifurcation points. Note that, if one computes less snapshots during the offline phase or use a worse discretization, it is very complex to obtain the diagram online, therefore 2365 is a reliable number in the considered scenario.}
\label{tab_time2}
\end{table}\\It is also important to remark that the obtained diagram is still very accurate. To prove it, we performed the same error analysis of the one-dimensional diagram obtaining the result in figure \ref{2d_error}. Such a figure has been obtained computing the relative error of the reduced solutions with respect to their full order counterpart obtained with the SEM. The average error is higher because following the online branches in a multi dimensional space is much more complex. However, the same behaviour shown in figure \ref{1d_error} is present, the error is characterized by two groups of peaks near the bifurcation points and by an oscillatory behaviour elsewhere. We decided to use only 3 values of $s$ to obtain a more clear visualization, such a choice does not influence the analysis since the most important parameter in diagram \ref{2d_diagr} is $\nu$.\\
Once more, we summarized in table \ref{tab2d} the average and maximum relative error according to the different regions of the diagram. The columns represent the same regions as in table \ref{tab1d} but the key parameter is the ratio between the viscosity and the scaling. Such a choice is justified by the fact that, with a fixed domain, this ratio determines the Reynolds number and, therefore, can be used to divide the different regions. It can be noted that, again, the error significantly increases near the bifurcation points, but that the entire diagram is, on average, still very accurate.
\begin{center}
\begin{table}[h]
\centering 
\begin{tabular}{c c c c c} 
\hline \\[-0.2cm]
Quantity&\hspace{0.15cm}$\frac \nu s\in\Omega_0^{\nu,s}$\hspace{0.15cm}&\hspace{0.15cm}$\frac  \nu s\in\Omega_1^{\nu,s}$\hspace{0.15cm}&\hspace{0.15cm}$\frac  \nu s\in\Omega_2^{\nu,s}$\hspace{0.15cm}&\hspace{0.15cm}$\frac  \nu s\in\Omega^{\nu,s}$\hspace{0.15cm}\\ [0.6ex]
\hline \\[-0.2cm]
Average error & $5.08\cdot10^{-4}$ & $1.02\cdot10^{-4}$ &$ 9.69\cdot10^{-4}$ &  $2.64\cdot10^{-4}$\\ 
Maximum error &$1.86\cdot10^{-2}$& $3.40\cdot10^{-3}$ & $9.49\cdot10^{-3}$ &  $1.86\cdot10^{-2}$\\[1ex] 
\hline\\[-0.8ex] 
\end{tabular}
\caption{Average and maximum relative error of the solutions of diagram \ref{2d_diagr} with respect to the corresponding full order solutions. The following notation has been used: $\Omega_0^{\nu,s}=[0.3,0.45]$, $\Omega_1^{\nu,s}=(0.45,0.9]\cup(1,1.25]$, $\Omega_2^{\nu,s}=(0.9,1]$, $\Omega^{\nu,s}=\Omega_0^{\nu,s}\cup\Omega_1^{\nu,s}\cup\Omega_2^{\nu,s}=[0.3,1.25]$.} 
\label{tab2d}
\end{table}
\end{center}\begin{figure}
\centering  
  \begin{subfigure}{9.5cm}
    \centering
    \includegraphics[width=9.5cm]{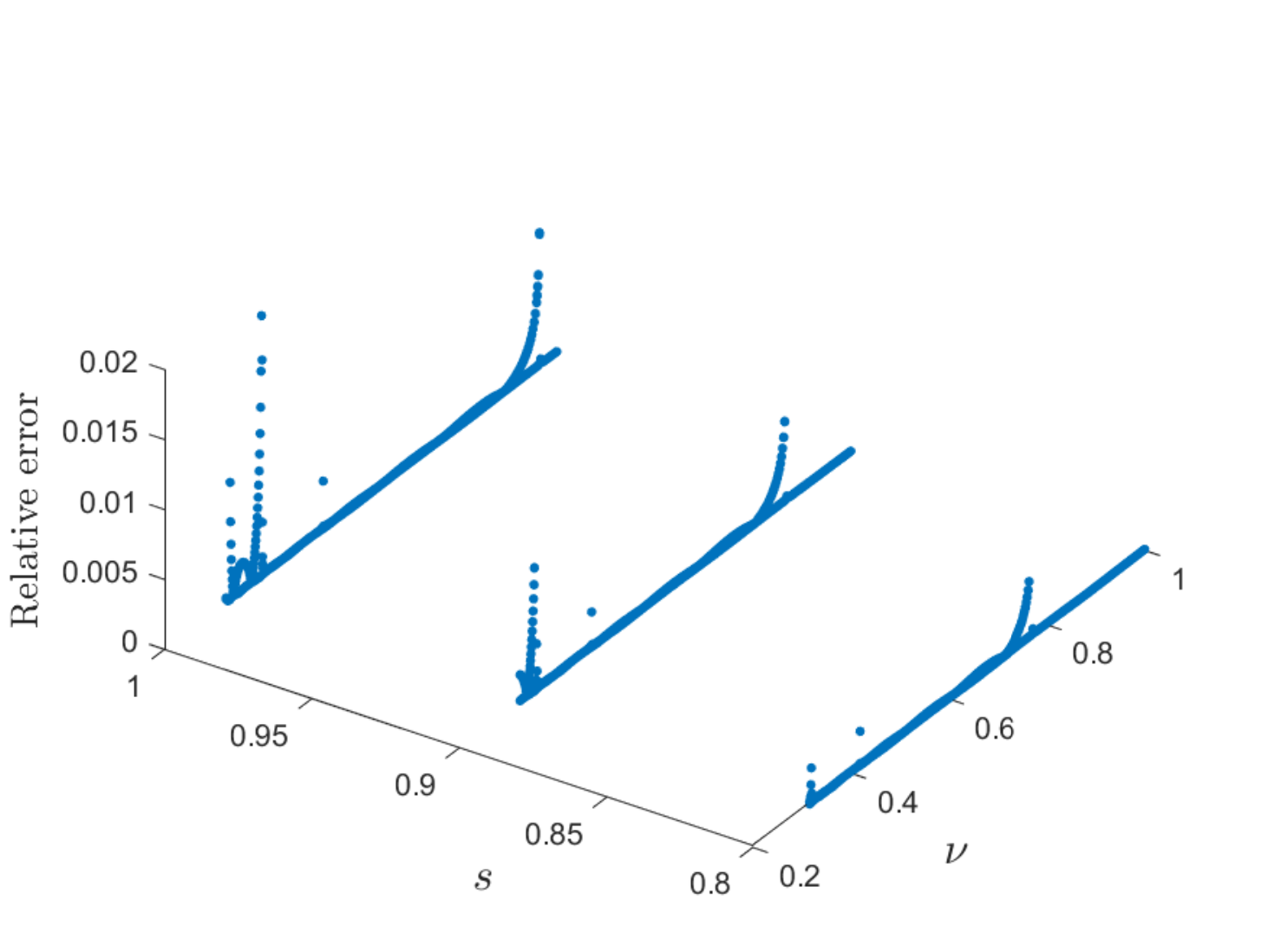}
    \caption{Uniform scale}
    \label{error_2d_unif}
  \end{subfigure}
  \begin{subfigure}{9.5cm}
    \centering
    \includegraphics[width=9.5cm]{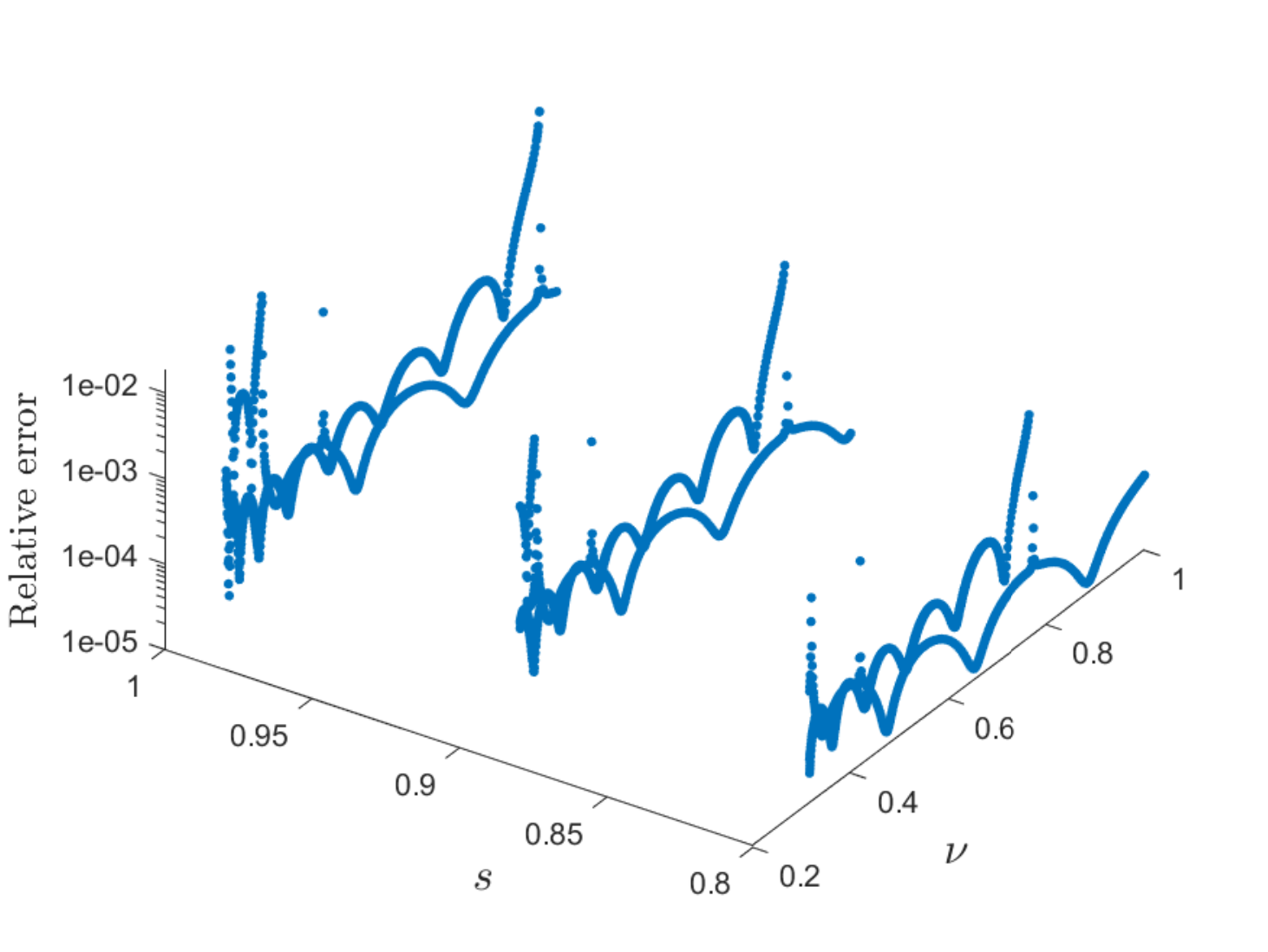}
    \caption{Logarithmic scale}
    \label{error_2d_log}
  \end{subfigure}
  \caption{Relative error with respect to full order solutions}
  \label{2d_error}
\end{figure}Finally, we are interested in understanding the effect of a geometrical variation on the bifurcation diagram. The corresponding parameter will be called $c_H$ and will represent a multiplicative factor of the inlet height used in figure \ref{mesh} and in the previous numerical results. When $c_H=1$, the domain used to compute the previous diagrams is considered and the inlet height is $2.5$. To be consistent with the velocity scaling factor $s$, $c_H$ will belong to the interval $[0.8,1]$. Therefore the inlet height varies in $[2,2.5]$ and the inlet Dirichlet boundary condition in equation \eqref{real_bc} has to be modified as follows:
\begin{equation*}
\uu  = \begin{bmatrix} u\\ v \end{bmatrix} = \begin{bmatrix} 20s(y_1^{c_H}-y)(y-y_0^{c_H}) \\ 0 \end{bmatrix} , \hspace{1cm} x=0,\text{ }y\in(y_0^{c_H},y_1^{c_H}),
\end{equation*} 
where $y_0^{c_H}$ and $y_1^{c_H}$ can be computed as:
\[
y_0^{c_H} = \dfrac{5+2.5}{2} - \dfrac{5-2.5}{2}h_C,\hspace{1cm}y_1^{c_H} = \dfrac{5+2.5}{2} + \dfrac{5-2.5}{2}h_C.
\]
As in the previous case, we computed, during the offline phase, a set of snapshots belonging to a small set of one-dimensional bifurcation diagrams. In particular, we computed them letting vary only $\nu$ and with $(s,c_H)\in\{0.8,1\}\times\{0.8,0.9,1\}$ (i.e. only 6 one-dimensional diagrams are obtained in the offline phase). Then, all the snapshots are used to construct a global reduced space by means of the POD and such a space is finally used to generate the corresponding three-dimensional diagram. In order to better visualize and compare the obtained diagram, we decided to compute six one-dimensional diagrams (in the $\nu$-direction) corresponding to parameter values $(s,c_H) \notin \{0.8,1\}\times\{0.8,0.9,1\}$ and show them together, with different colours, in figure \ref{3d_diagr}.\\
In particular, the red curves are obtained with $c_H=0.95$ and the blue ones with $c_H=0.85$. It can be observed that the bifurcation diagrams are similar but that the critical points are anticipated as $c_H$ increases. Once more, we analyzed the accuracy of the involved solutions with respect to the corresponding full order ones. The obtained error diagram is shown in figure \ref{3d_error}; it can be noted that it is consistent with the ones in figures \ref{1d_error} and \ref{2d_error}.
\begin{figure}
\centering
\includegraphics[height=11.8cm, width=11.8cm, angle=0, keepaspectratio]{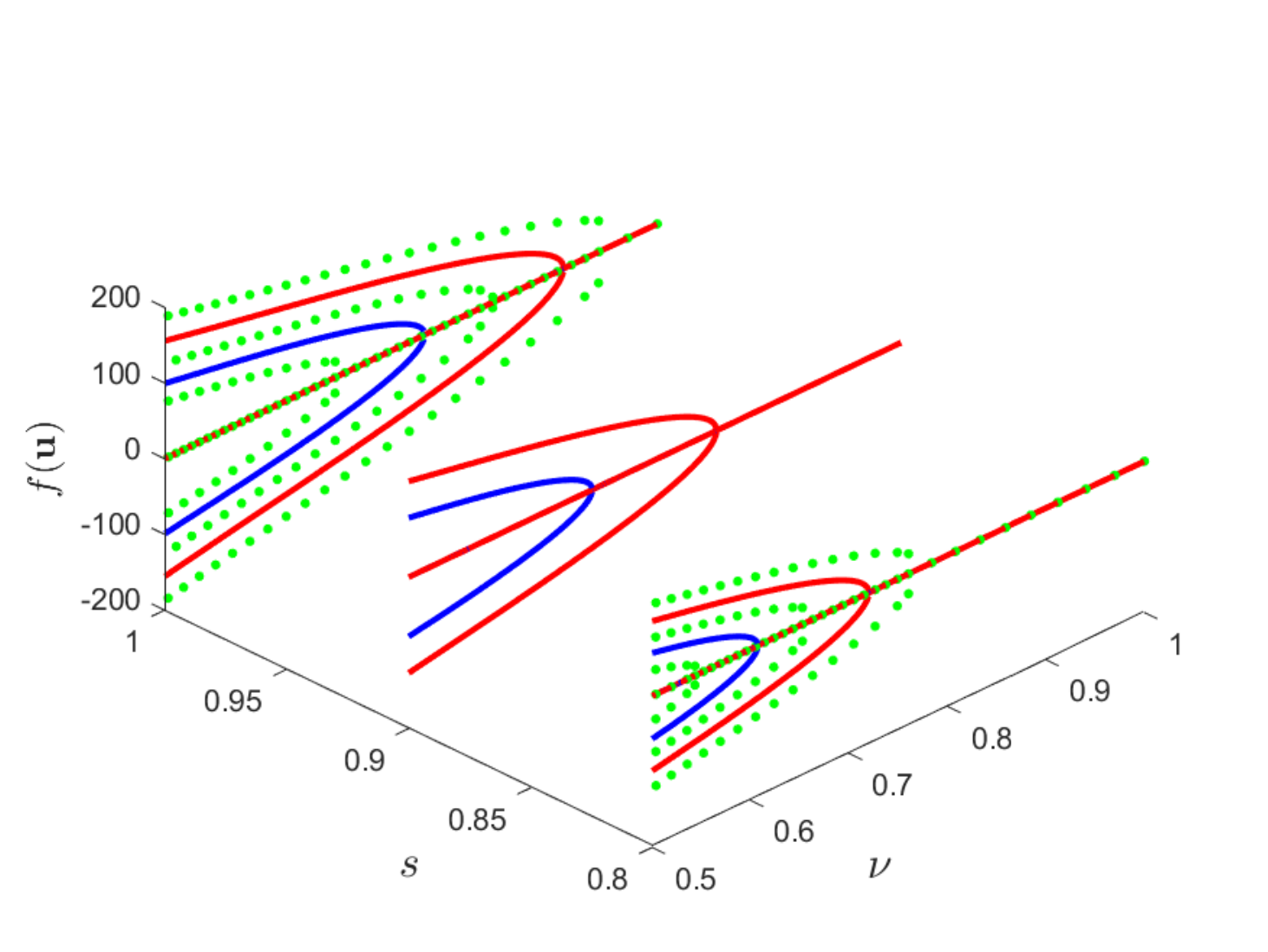}
\caption{Bifurcation diagram with three parameters efficiently obtained in the online phase. The green dots represent the snapshots, while the red and blue curves the bifurcation diagrams obtained with $c_H=0.95$ and $c_H=0.85$, respectively}
\label{3d_diagr}
\end{figure}
\begin{figure}
\centering
\includegraphics[height=11.8cm, width=11.8cm, angle=0, keepaspectratio]{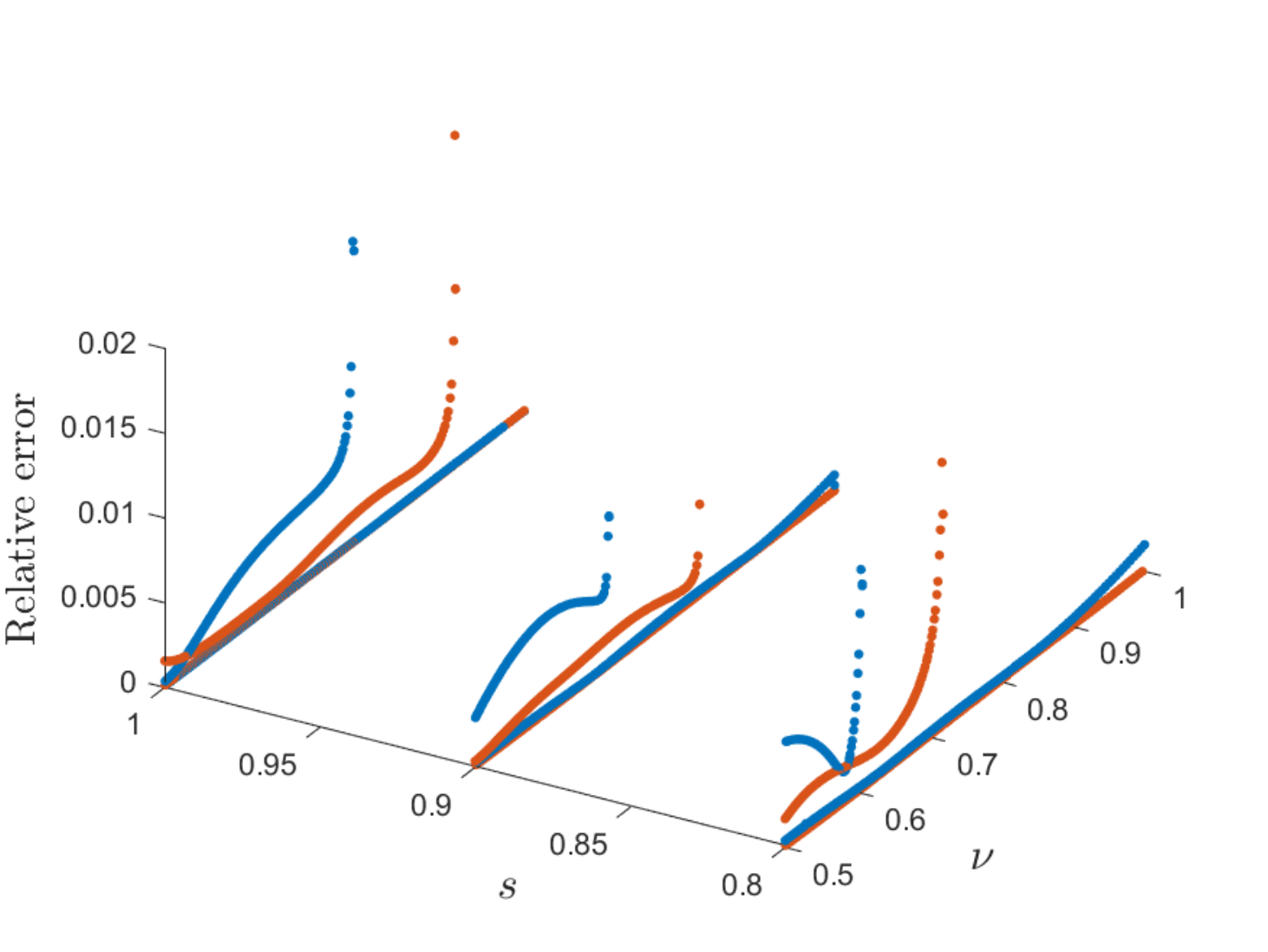}
\caption{Relative error with respect to full order solutions }
\label{3d_error}
\end{figure}

\section{Conclusions and perspectives}
\label{concl_chapter}
In this work we described an approach to efficiently compute bifurcation diagrams with more than one parameter. The diagrams are discretized exploiting an elaborated deflated continuation method characterized by two versions of the continuation method and a deflation one which associated steps can be adaptively modified to improve its efficiency and effectiveness. Moreover, we decided to use the reduced basis method to drastically reduce the computational cost of the process and, in order to increase the efficiency of its offline phase, we adopted the spectral element method.\\
The main advantages of the described method are the following ones:
\begin{itemize}
  \item the bifurcation diagram is automatically computed, a priori knowledge of the phenomenon is not required;
  \item when more than one parameter is involved, the computational cost of the entire process is drastically reduced thanks to the reduced basis method. In fact it allows one to compute any solution at a cost independent of $N^\delta$ and to compute offline only few full-order solutions;
  \item the offline phase can be further improved exploiting high-order methods as the SEM, because they require a lower number of degrees of freedom than low-order methods and, since the supports of the bases are wider, they are more similar to the ones used in the online phase;
  \item with the continuation method and the deflation one it is possible to compute bifurcation diagrams with an arbitrary number of bifurcation points, whose nature is not known a priori.
\end{itemize}
However, it is important to recognize that the stability of the online solver should be improved before applying such a technique in real scenarios. In fact, in order to obtain the convergence during the online phase, we had to increase the order of the polynomials in the offline one and to properly select some parameters. A more sophisticated approach could involve the supremizers \cite{supremizer} to directly improve the stability or the localized reduced basis method \cite{local_rom}, that can reduce the noise in each basis. Moreover, we expect that such an issue can worsen when other parameters, mainly the geometrical ones, are involved and, therefore, the described method has to be tested on three-dimensional and realistic geometries before being applied in real scenarios. It should be noted that the fact that the SEM relies on coarse meshes is not restrictive, in fact a complex geometry can be accurately approximated by means of curved elements \cite{martin_curve_walls}. Moreover, to further increase the efficiency and the stability of the proposed technique, one can generalize it with the reduced basis element method \cite{rbem}. In such a method the wide supports of the SEM basis functions are exploited to generate a reduced space associated with each element or with groups of elements. We remark that combining complex three-dimensional simulations as the ones in \cite{cardio} with the described advanced numerical techniques is still very challenging. Therefore the proposed method cannot be immediately used to understand the flow features and detect whether there are anomalies related to the blood flow.\\
Finally, we have only studied steady bifurcations where a finite number of solutions existed, however, more advanced studies and tools are required to capture unsteady bifurcations (as the Hopf bifurcations \cite{bif_t}, obtained, for instance, in \cite{sobey1986bifurcations} for slightly higher Reynolds numbers in a similar geometry) or phenomena characterized by infinite solutions. For instance, in \cite{deflation_groups}, a different deflation operator has been implemented to deflate entire groups of solutions, characterized by their symmetry, at the same time. Moreover, we highlight that if one is interested in the nature of the bifurcation points, it is convenient to analyze the behaviours of the eigenvalues that lead to the bifurcations (see \cite{pichi_schro}).

\section*{Acknowledgements}
We acknowledge the support by European Union Funding for Research and Innovation - Horizon 2020 Program - in the framework of European Research Council Executive Agency: Consolidator Grant H2020 ERC CoG 2015 AROMA-CFD project 681447 ``Advanced Reduced Order Methods with Applications in Computational Fluid Dynamics"  (P.I. Prof. Gianluigi Rozza).\\ 
We also acknowledge the INDAM-GNCS project ``Advanced intrusive and non-intrusive model order reduction techniques and applications", 2019.

\bibliographystyle{spmpsci}      
\bibliography{bibliography}   

\end{document}